\documentclass[10pt]{article}

\usepackage{amssymb,latexsym,amsmath,makeidx}
\newtheorem{theo}{Theorem}[section]
\newtheorem{lemm}[theo]{Lemma}

\newtheorem{prop}[theo]{Proposition}
\newtheorem{defi}[theo]{Definition}

\numberwithin{equation}{section}
\newenvironment{Proof}{\removelastskip\vskip12pt plus
1pt\noindent\em Proof. \rm}{\hspace*{\fill}\ $\Box$\vskip12pt
plus 1pt}
\newcommand{\End}{{\rm End}}

\newcommand{\Der}{{\rm Der\,}}
\newcommand{\ad}{{\rm ad\,}}
\newcommand{\rad}{{\rm rad\,}}
\newcommand{\gr}{{\rm gr \,}}
\newcommand{\N}{{\Bbb N}}

\newcommand{\F}{{\Bbb F}}
\newcommand{\s}{\frak{sl}}

\begin{document}

\title{\bf Lie algebras of small dimension}

\author{\bf H. Strade}
\date{}
\maketitle

\renewcommand{\thefootnote}{\fnsymbol{footnote}}
\footnotetext[1]{
{\it Revision}). Primary  17B55, 17B65}

\begin{abstract}
We present a list of all isomorphism classes of nonsolvable Lie algebras of dimension $\le 6$ over a finite field.
\end{abstract}

\section{Introduction}

During the last years people have been very successful in classifying groups up to an order of 2000 and compile them into the computer program GAP and MAGMA. This makes it possible for group theorists to test conjectures or find counterexamples, and this has turned out to be a very fruitful procedure.

It seems therefore a desirable task to do something similar for Lie algebras. There is an old result of Zassenhaus and Patera (\cite{ZP}), who classify  solvable Lie algebras up to dimension 6 over finite fields.  Quite recently, W. de Graaf corrected and extended their work  (\cite{G}). In particular, he listed all at most 4-dimensional  solvable Lie algebras over arbitrary fields. C. Schneider worked on nilpotent Lie algebras over small fields (\cite{Sch}). He classified all such algebras of dimension at most 6 over finite fields of characteristic $>2$. Due to C. Schneider, W. de Graaf and others a computer program for nilpotent Lie algebras is under way, which at present might allow an inductive classification for probably at most 6-dimensional nilpotent Lie algebras over finite fields.

Appearantly, not much is known so far about nonsolvable Lie algebras over finite fields. It turns out that for these classes of Lie algebras methods are useful which have been developed in the course of the classification of the finite dimensional simple Lie algebras over algebraically closed fields. In this paper we classify all nonsolvable Lie algebras over finite fields up to dimension 6. 

\bigskip

We remind the reader to some concepts and results from Lie algebra theory over fields of positive characteristic $p$. The reader may find details and proofs of the following facts in  \cite{SF} and \cite{St04}.

A $p$-{\em envelope} $(G,[p],\iota)$ of a Lie algebra $H$ is a triple consisting of a restricted Lie algebra $(G,[p])$ and an injective homomorphism $\iota:H\hookrightarrow G$ such that the algebra generated by $\iota(H)$ and $[p]$ is $G$. Such $p$-envelopes always exist, and finite dimensional $p$-envelopes exist if $H$ is finite dimensional. The $p$-envelope of a subalgebra $H$ of a restricted Lie algebra $(G,[p])$ is the algebra $H_{[p]}\subset G$ generated by $H$ and $[p]$ (with $\iota$ the inclusion).
\\
Let $H\subset G$ be finite dimensional Lie algebras, $(G_p,[p],\iota)$ be a $p$-envelope of $G$ and $H_{[p]}$ the $p$-envelope of $\iota(H)$ in $G_{[p]}$. Then $\dim\,H_{[p]}/H_{[p]}\cap C(G_{[p]})<\infty$ (as $G$ is finite dimensional), and  $H_{[p]}/H_{[p]}\cap C(G_{[p]})$ carries a natural $p$-mapping inherited by the $p$-mapping of $H_{[p]}$. The maximal dimension of tori in $H_{[p]}/H_{[p]}\cap C(G_{[p]})$ is called the {\em toral rank} $TR(H,G)$ of $H$ in $G$. This number is independent of the $p$-envelope chosen. If $H=G$, then $TR(G):=TR(G,G)$ is called the {\em absolute toral rank} of $G$. This concept is of importance in the classification theory of finite dimensional simple Lie algebras over fields of positive characteristic because in general finite dimensional restricted Lie algebras contain maximal tori of various dimensions. In contrast, if $G$ is nilpotent then any finite dimensional $p$-envelope $G_{[p]}$ is also nilpotent and contains a unique maximal torus.

Let $N$ be a nilpotent subalgebra of a finite dimensional Lie algebra $G$. The {\em primary decomposition} of $G$ with respect to $N$ is described in \cite[\S 1.4]{SF} (in particular in Theorem 1.4.3). If the ground field $k$ is algebraically closed, irreducible polynomials over $k$ are linear and therefore there are eigenvalue functions $\alpha:N\to k$ such that
$$G=\oplus_\alpha G_\alpha(N),\quad G_\alpha(N):=\{g\in G\mid (\ad x-\alpha(x){\rm Id})^{\dim G}(g)=0\ \forall x\in N\}.$$
These {\em eigenvalue functions} are not necessarily linear functions on $N$. Sometimes they are called {\em extended roots}. The set of extended roots (including 0!) is denoted by $\Gamma(G,N)$.
\\
If the ground field $k$ is arbitrary, then we don't necessarily have eigenvalue functions available. Nevertheless, the {\em Fitting decomposition} of $G$ with respect to an element $x$ is of major importance. There is always the decomposition
$$G=G_0(\ad x)\oplus G_1(\ad x)$$
where $x$ acts nilpotently on $G_0(\ad x)$ and invertibly on $G_1(\ad x)$.
\\[2mm]

Let ${\cal O}(m)$ denote the commutative and associative $k$-algebra with unit element defined by generators $$x_i^{(r)},\quad 1\le i\le m,\ r\ge0$$ and relations 
$$x_i^{(0)}=1,\quad x_i^{(r)}x_i^{(s)}={r+s\choose r}x_i^{(r+s)}.$$
Put $x_i:=x_i^{(1)}$, $x^{(a)}:=x_1^{(a_1)}\cdots x_m^{(a_m)}$ for $a\in\N^m$. Then $(x^{(a)}\mid a_i\ge0)$ is a basis of ${\cal O}(m)$. For any $m$-tuple $\underline n=(n_1,\ldots,n_m)\in\N^m$ we set 
$${\cal O}(m;\underline n):=span\{x^{(a)}\mid 0\le a_i\le p^{n_i}-1\},$$
and observe that this is a subalgebra of ${\cal O}(m)$ of dimension $p^{n_1+\cdots+n_m}$. 
\\[2mm]
The following subset of derivations
$$W(m):=\{D\in \Der{\cal O}(m)\mid D(x_i^{(r)})=x_i^{(r-1)}\ \forall i=1,\ldots,m,\ \forall r>0\}$$
is a Lie algebra, the algebra of {\em special derivations}. Accordingly, 
$$W(m;\underline n):=W(m)\cap \Der{\cal O}(m;\underline n)$$
is a Lie algebra of {\em Witt type}. It has finite dimension $mp^{n_1+\cdots+n_m}$. The {\em partial derivatives} are defined by $\partial_i(x_j^{(r)})=\delta_{i,j}x_j^{(r-1)}$ for arbitrary $i,j=1,\ldots,m$ and all $r>0$.
\\
These definitions can be done over arbitrary fields $k$. Sometimes in our context the ground field will play an important role. In such a situation we will insert the notion of $k$ by writing ${\cal O}(k| m;\underline n)$ and $W(k| m;\underline n)$. We will, however, suppress the notion of the ground field whenever it is clear which ground field is meant.
\\
For the definition of the series of Cartan type Lie algebras we refer to the Definitions 4.2.1 and 4.2.4 of \cite{St04}. They only occur once in the present paper. It is sufficient for our purposes to mention that $W(1;\underline 1)$ is the  only Lie algebra of Cartan type over a field of characteristic $p\ge3$, which has dimension $\le6$. Also, $W(m;\underline n)$ is simple except if $m=1$ and the ground field has characteristic $p=2$. In the latter case $W(1;\underline n)^{(1)}$ is simple of dimension $2^n-1$.

\bigskip

The precedure of our investigations is as follows. We first determine the semisimple Lie algebras of dimension $\le 6$ over an algebraically closed field. Then we describe the forms of these algebras, i.e., Lie algebras $L$ over a finite field $k$ for which $L\otimes _k\bar k$ ($\bar k$ the algebraic closure of $k$) is isomorphic to one of these semisimple Lie algebras. In a third step we consider the extensions 
$$0\to \rad(L)\to L\to L/\rad(L)\to 0.$$
The results of the first and second step can be formulated with only little dependency  on the characteristic of the ground field and the dimension of the Lie algebra. Up to dimension 5 the number of isomorphism types of arbitrary Lie algebras is strongly limited. However, there is a rather long list of isomorphism types of 6-dimensional Lie algebras. In particular,  the characteristic 3 case is somewhat massy.

\section{Semisimple Lie algebras}

We mention an important principle to construct ideals in a Lie algebra.

\begin{lemm}(cf. \cite[Proposition 1.3.5]{St04})\label{2.1} Let $G$ be an arbitrary Lie algebra which is finite dimensional over an algebraically closed field of positive characteristic, and $N\subset \Der G$ a nilpotent subalgebra. Then 
$$\sum_{\mu\in\Gamma(L,N)\setminus\{0\}}\,G_\mu+\sum_{\lambda,\mu\in\Gamma(L,N)\setminus\{0\}}\,[G_{\lambda},G_\mu]$$
is an ideal of $G$.
\end{lemm}
\begin{Proof}
Note that the Fitting-1-component of $G$ with respect to $N$ is $G_1(N)=\sum_{\mu\in\Gamma(L,N)\setminus\{0\}}\,G_{\mu}$. The statement is a direct consequence of this observation and Proposition 1.3.5 of \cite{St04}.
\end{Proof}

 This lemma will often be applied  when $N=T$ is a torus of $\Der G$ or $N$ is a subtorus $T_0:=\cap_{\mu\in\Gamma_0\subset\Gamma(G,T)}\,\ker\mu $ of a given torus $T$.

\vspace{3mm}

Lie algebras of dimension $<3$ as well as 3-dimensional Lie algebras containing a proper ideal are solvable. Therefore  every at most 3-dimensional nonsolvable Lie algebra is  simple.

\begin{theo}\label{2.2}
Let $G$ be a simple Lie algebra over an algebraically closed field $F$ of characteristic $p>0$. Suppose $\dim G\le6$. Then $G$ is 3-dimensinonal or $p=5$ and $G\cong W(1;\underline 1)$.
\end{theo}
\begin{Proof}
Let $H$ be a Cartan subalgebra of $G$, let $H_{[p]}$ be the $p$-envelope of $H$ in $\Der G$, and $T$ denote the maximal torus of the nilpotent algebra $H_{[p]}$. Let $G=\sum_{\mu\in\Gamma(G,T)}G_\mu$ denote the root space decomposition of $G$ with respect to $T$.
\\[2mm]
(1) Consider the case $p>2$.\\[2mm]
(i) Suppose $\dim T>1$. The simplicity of $G$ implies
$H=\sum_{\mu\ne0}[G_{-\mu},G_\mu]$ (Lemma \ref{2.1}). 
As $H\ne\{0\}$, there is a nonzero  root  $\alpha$ for which $-\alpha$ is a  root. The present assumption $\dim T>1$ yields that $T\cap\ker \alpha\ne\{0\}$. Choose a toral element $t\in T\cap\ker\alpha$ and decompose 
$G=\sum_{i\in \F_p}G_i$
into the $t$-eigenspaces. Note that $H+G_{-\alpha}+G_\alpha\subset G_0$, whence $\dim G_0\ge3$. The simplicity of $G$ yields $G_0=\sum_{i\ne0}[G_{-i},G_i]$, and therefore 
$$\sum_{i=1}^{(p-1)/2}(\dim G_{-i})(\dim G_i)\ge \dim G_0\ge3.$$
On the other hand, $\sum_{i=1}^{p-1}\dim G_i=\dim G - \dim G_0\le3$,  which contradicts the former inequality. \\[2mm]
(ii) As a result,  $\dim T=1$ holds. This means in other words, that $H$ has toral rank $TR(H,G)=1$ in $G$. Now \cite[Theorem\ 9.2.11]{St04} shows that $G$ is one of 
$$\s(2),\ W(1;\underline n),\ H(2;\underline n;\Phi)^{(2)}.$$
Note that $\dim H(2;\underline n;\Phi)^{(2)}\ge p^2-2\ge7$
and 
$\dim W(1;\underline n)=p^n$.
Consequently, no Hamiltonian algebra occurs, and the only algebra of Witt type occurring is $W(1;\underline 1)$ for $p\le5$ (which is $p$-dimensional).
\\[2mm]
(2) Consider the case $p=2$.\\[2mm]
(i) Suppose $\dim T>1$. The simplicity of $G$ implies
$$H=\sum_{\mu\ne0}[G_{\mu},G_\mu].$$
Therefore there is a nonzero root $\alpha$ for which $\dim G_\alpha\ge2$. Choose a toral element $t\in T\cap\ker\alpha$ and decompose 
$G=G_1\oplus G_0$
into the $t$-eigenspaces. As before we obtain $\dim G_0\ge3$, $\dim G_1\le3$, and 
$$\frac{1}{2}(\dim G_{1})(\dim G_1-1)=\dim\,G_1\wedge G_1\ge \dim G_0.$$
But then 
$$\dim G_1=3,\ \dim G_\alpha=2,\ \dim H=1.$$
Note that $G_0$ and $G_1$ are $T$-invariant. The simplicity of $G$ implies that $G=G_1+[G_1,G_1]$ (Lemma \ref{2.1}). Then $[G_1,G_1]=G_0$ and one gets $H\subset [G_1,G_1]$. Therefore there is a root $\beta$ with  $[\dim G_\beta\cap G_1,\dim G_\beta\cap G_1]\ne\{0\}$ and hence $\dim G_\beta\cap G_1\ge2$. The simplicity of $G$ also implies that   there is a root $\gamma$  for which $[G_\gamma\cap G_1,G_\alpha]\ne\{0\}$ (otherwise $[G_\alpha,G_1]=\{0\}$ and  the equality $G=G_1+[G_1,G_1]$ would imply that $G_\alpha$ is contained in the center of $G$). Hence $\alpha+\gamma$ is a root on $G_1$. For dimension reasons the roots $\beta$, $\gamma$, $\alpha+\gamma$ cannot be distinct. Therefore it is only possible that $\beta=\gamma$ or $\beta=\alpha+\gamma$.
Again by dimension reasons we conclude that  $\Gamma(G,T)=\{0,\alpha,\beta,\gamma\}$. Lemma \ref{2.1} also implies that for every nonzero root $\kappa$ one has $G_\kappa=\sum_{\lambda\ne\kappa}\,[G_\lambda,G_{\kappa+\lambda}]$.
 Since $\dim G_1=3$, we now have
\begin{align*}&G=H\oplus G_\alpha\oplus G_\beta\oplus G_{\alpha+\beta},\\
&\dim G_\alpha=\dim G_\beta=2,\ \dim G_{\alpha+\beta}=1,\\
&[G_\alpha,G_\alpha]=[G_\beta,G_\beta]=H=:Fh,\\&[G_\alpha,G_{\alpha+\beta}]=G_\beta,\ [G_\beta,G_{\alpha+\beta}]=G_\alpha,\ [G_\alpha,G_\beta]=G_{\alpha+\beta}.\end{align*}
Choose a basis $e_\alpha,f_\alpha$ of $G_\alpha$ such that $[h,e_\alpha]=\alpha(h)e_\alpha$. Adjusting $f_\alpha$ by a nonzero scalar we may assume that
$ [e_\alpha,f_\alpha]=h.$
Since $\dim\, G_{\alpha}=\dim\, G_{\beta}=2$ and $\dim\, [G_{\alpha},G_{\beta}]=1$ there are linearly independent elements $e_\beta,f_\beta\in G_\beta$ for which
$$[e_\alpha,e_\beta]=[f_\alpha,f_\beta]=0.$$
As $e_\alpha$ is an $h$-eigenvector we may take $e_\beta$ as an $h$-eigenvector as well. Adjusting $f_\beta$ by a scalar we obtain in total
\begin{align*}
&[h,e_\alpha]=\alpha(h)e_\alpha,&& [e_\alpha,f_\alpha]=h,\\& [h,e_\beta]=\beta(h)e_\beta,&& [e_\beta,f_\beta]=h,\\& [e_\alpha,e_\beta]=[f_\alpha,f_\beta]=0.\end{align*}
Since $ [e_\alpha,f_\beta],[f_\alpha,e_\beta]\in G_{\alpha+\beta}$ and $G_{\alpha+\beta}$ is 1-dimensional, 
\begin{align*}
0&=[[e_\alpha,f_\beta],[f_\alpha,e_\beta]]=[[e_\alpha,[f_\alpha,e_\beta]],f_\beta]+[e_\alpha,[f_\beta,[f_\alpha,e_\beta]]]\\
&=[[[e_\alpha,f_\alpha],e_\beta],f_\beta]+[e_\alpha,[f_\alpha,[f_\beta,e_\beta]]]\\&=[[h,e_\beta],f_\beta]+[e_\alpha,[f_\alpha,h]]=[[h,e_\beta],f_\beta]+[f_\alpha,[e_\alpha,h]]\\&=(\beta(h)+\alpha(h))h.
\end{align*}
But then $\alpha+\beta$ vanishes on $H$. Then it vanishes on $T$, which is impossible.
 \\[2mm]
 (ii) As a result, $\dim T=1$. Then there is only one nonzero root $\alpha$,
 $$G=H\oplus G_\alpha.$$
 Suppose $H^{(1)}$ acts nonnilpotently on $G_\alpha$. Then $\dim H\ge3$, and there is a  composition factor of a $H$-composition series of $G_\alpha$ having dimension $>1$. Every irreducible $H$-module has  2-power dimension (\cite[Corollary 3.2.8]{St04}). As $\dim G_\alpha=\dim G-\dim H\le3$, it can only be that $\dim G_\alpha=2$. But then $H=[G_\alpha,G_\alpha]$ would be 1-dimensional. As a consequence, $H^{(1)}$ acts nilpotently on $G_\alpha$.
 \\[2mm]
(iii) Let $Q$ denote a maximal subalgebra of $G$ containing $H$. Then $L=G_\alpha+Q$ and  $Q=H+G_\alpha\cap Q$.  Since $[G_\alpha,G_\alpha]\subset H\subset Q$ and $H^{(1)}$ acts nilpotently on $G$, one has that $Q^{(1)}\subset H^{(1)}+G_\alpha\cap Q+[G_\alpha\cap Q,G_\alpha\cap Q]$ acts nilpotently on $G/Q$.  Then $\{x\in G\mid [x,Q^{(1)}]\subset Q\}\supsetneq Q$ and is $Q$-invariant, and therefore $G/Q$ contains a $Q$-eigenvector,
$$\exists y\in G\setminus Q,\quad [Q,y]\subset Fy+Q.$$
We obtain that $Q+Fy$ is a subalgebra of $G$, and the maximality of $Q$ implies $G=Q+Fy$. In particular, $G/Q$ is $Q$-irreducible.
 Define a standard filtration of $G$ (cf. \cite[p. 168]{St04}),
 $$G_{(-1)}:=G,\ G_{(0)}:=Q,\ G_{(i+1)}:=\{x\in G_{(i)}\mid [G,x]\subset G_{(i)}\}\ \text{for }i\ge0.$$
 Note that $\dim G/G_{(0)}=1$, and therefore $\dim G_{(i)}/G_{(i+1)}=1$ holds whenever $G_{(i)}$ is nonzero. Set $q:=\dim G-2$ and let $\gr G:=\oplus_{i=-1}^qG_{(i)}/G_{(i+1)}$ denote the associated graded algebra. Choose a homogeneous basis for $\gr G$,
 $$\gr_{-1}G=:Fe_{-1},\ \gr_q G=:Fe_q,\ e_{i}:=[e_{-1},e_{i+1}]\ \text{for }0\le i\le q-1.$$
 We may adjust $e_q$ by a scalar so that $[e_0,e_{-1}]=e_{-1}$. Then
 \begin{align*}
\gr_i G &=Fe_i,&&i=-1,\ldots,q,\\ [e_{-1},e_i] &=e_{i-1},&&i=0,\ldots,q,
\\  [e_i,e_j] &\in \gr_{i+j} G,&&i+j\ge-1,
\end{align*}
and one easily proves by induction that $[e_0,e_i] =ie_i$ holds for all $i=-1,\ldots,q$.
The simplicity of $G$ implies $\dim G\ge3$, whence $q\ge1$. The case $q=1$ is the claim. 
 Suppose $q\ge2$. Then 
 $$[e_{-1},[e_1,e_2]]=[[e_{-1},e_1],e_2]]+[e_{1},[e_{-1},e_2]]=[e_0,e_2]+[e_1,e_1]=0,$$
 whence 
 $$[e_1,e_2]=0.$$
 If, moreover, $q>2$, then 
 $$[e_{-1},[e_1,e_3]]=[e_0,e_3]+[e_1,e_2]=e_3,$$
 whence $[e_1,e_3]\ne0$ and $q\ge4$. Next,
 $$[e_{-1},[e_1,e_4]]=[e_0,e_4]+[e_1,e_3]=[e_1,e_3]\ne0.$$
 Then $[e_1,e_4]\ne0$ and $q\ge5$. But then $\dim G\ge 7$, a contradiction. \\[2mm]
 Consequently, $q=2$ under the present assumption and $\dim G=4$. If $\dim\, H$ would be bigger than 1, then only $\dim\, G_\alpha=2$ is possible and this gives $\dim\, H\le \dim\, [G_\alpha,G_\alpha]\le1$. This contradiction shows that $H=:Fh$ is 1-dimensional.
Note that $h\in G_{(0)}\setminus G_{(1)}$. Hence $F(h+G_{(1)})=Fe_0$. Choose a preimage $h'$ of $e_2$ in $G_{(2)}$. As $[e_0,e_2]=0$, one obtains $[h,h']\in G_{(3)}=\{0\}$. Then $h'\in H$. But $H$ is 1-dimensional, and this contradiction completes the proof of the theorem.
\end{Proof}

Using this result we can determine the semisimple Lie algebras over a finite field to some extent.

\begin{defi}\label{2.3}
Let $L$ be an arbitrary Lie algebra and $M$ be an irreducible $L$-module. Set $${\cal C}(M,L):=\{f\in\End M\mid f(x.m)=x.f(m)\ \forall x\in L,\ \forall m\in M\}.$$
 If $L$ is simple, then set ${\cal C}(L):={\cal C}(L,L)$.
\end{defi}

Schur's lemma shows that ${\cal C}(M,L)$ is a division algebra. Suppose the ground field $k$ is finite and $M$ is finite dimensional over $k$. As ${\cal C}(M,L)$ is finite dimensional over the finite field $k$, it is finite. Now Wedderburn's theorem proves that ${\cal C}(M,L)$ is a finite field extension of $k$. Note that $M$ is a vector space over ${\cal C}(M,L)$.

\begin{theo}\label{2.4} Let $L$ be a semisimple Lie algebra over a finite field $k$ of dimension $\le 6$. Let $\bar k$ denote the algebraic closure of $k$. \\[2mm]
If $p=2$, then one of the following holds.
\begin{enumerate}
\item $L=L_1\oplus L_2$ is the direct sum of simple 3-dimensional Lie algebras $L_1,L_2$;
\item $L$ is simple,  ${\cal C}(L)/k$ is a field extension of degree 2, and $L$ is 3-dimensional over ${\cal C}(L)$;
\item $L$ has a simple 3-dimensional ideal $S$, and $L\subset \Der S$.
\end{enumerate}
If $p\ge3$, then one of the following holds.
\begin{enumerate}
\item $L=L_1\oplus L_2$, where $L_i$ are simple and $L_i\otimes_k\bar k\cong \s(2,\bar k)$ ($i=1,2$);
\item $L$ is simple,  ${\cal C}(L)/k$ is a field extension of degree 2, and $L\otimes_{{\cal C}(L)}\bar k\cong \s(2,\bar k)$;
\item  $L\otimes_k\bar k\cong\s(2,\bar k)$;
\item $p=5$ and $L\cong W(k|1;\underline 1)$.
\end{enumerate}

\end{theo}
\begin{Proof}
(a) Let $\oplus_{i=1}^t\,L_i$ denote the socle of $L$, i.e., the direct sum of all minimal ideals. The minimality of $L_i$ implies that every $L_i$ is $L$-simple. Definition \ref{2.3} applies for $M=L_i$. Then $k_i:={\cal C}(L_i,L)$ is a finite field extension of $k$. Observe that $L_i$ is a nonsolvable Lie algebra over $k_i$. Therefore
$$6\ge \dim_k L\ge \sum_{i=1}^t \dim_kL_i=\sum_{i=1}^t (\dim_{k_i} L_i)(\dim_k k_i)\ge 3\sum_{i=1}^t (\dim_k k_i).$$
Consequently, if $t\ge 2$, then  $L=L_1\oplus L_2$, $k_i=k$, and $\dim_k L_i=3$   for $i=1,2$. In particular, $L_1$, $L_2$ are simple Lie algebras in their own right.
\\
 If $t=1$, then only $k_1=k$, or $$\dim_k k_1=2\quad \text{and}\quad  L=L_1\quad \text{and}\quad \dim_{k_1}L=3$$ are possible. In the latter case $L$ is a Lie algebra over $k_1$, hence is simple (as it is 3-dimensional nonsolvable).
\\[2mm]
(b)  If $t=2$, then (a)  implies that $L_1,L_2$ are central simple over $k$. It is well known that in this case $L_i\otimes_k\bar k$ again is simple. Then $L\otimes_{k}\bar k\cong (L_1\otimes_{k}\bar k)\oplus (L_2\otimes_{k}\bar k)$ is the direct sum of simple Lie algebras of dimension 3 over $\bar k$. If $p\ge3$, this gives $ L_i\otimes_{k}\bar k\cong \s(2,\bar k)$. \\
Similarly, if $t=1$ and $k_1\ne k$, then  $L\otimes_{k_1}\bar k$ is simple and 3-dimensional. As above, if $p\ge 3$ then this algebra is isomorphic to $\s(2,\bar k)$.\\[2mm]
(c) We finally consider the case that $L$ has a unique minimal ideal $L_1=:S$, and $k_1=k$. As $S$ is $L$-simple, the associative algebra generated by ${\rm ad}_SL$ is $\End_kS$ (Wedderburns's theorem, observe that $\End_L(S)={\cal C}(S,L)=k$). Then the associative algebra generated by ${\rm ad}_{(S\otimes_k\bar k)}(L\otimes_k\bar k)$ is $\End_{\bar k}(S\otimes_k\bar k)$, and therefore $S\otimes_k\bar k$ is $(L\otimes_k\bar k)$-simple. Due to Block's result (see \cite[Corollary 3.3.3]{St04}), $S\otimes_k\bar k\cong S'\otimes_{\bar k} {\cal O}_{\bar k}(m;\underline n)$, where $S'$ is a simple Lie algebra over $\bar k$ and $m\ge0$ and $\underline n\in \N^m$. For dimension reasons it is only possible that $m=0$ and $S'=S\otimes_k\bar k$, or $m=1$ and $\underline n=\underline 1$ and $p=2$ and $S'$ is 3-dimensional. However,  in the latter case  a dimension argument yields $L\otimes_k\bar k\cong S'\otimes {\cal O}_{\bar k}(1;\underline 1)\cong S\otimes_k\bar k$, whence $L=S$ and $S$ is central simple over $k$. But then $S\otimes_k\bar k$ would be simple, a contradiction. \\
As a result, $m=0$ and $S'=S\otimes_k\bar k$ is a simple Lie algebra over $\bar k$.   Theorem \ref{2.2} shows that $S\otimes_k\bar k$ is 3-dimensional or is isomorphic to $W(\bar k|1;\underline 1)$ (in case $p=5$). In the first case $S$ is 3-dimensional, and since $L$ is semisimple it embeds into $\Der S$. This solves the case $p=2$. 

Now assume $p\ge3$.
As every 3-dimensional Lie algebra over an algebraically closed field of characteristic $p\ge3$ is isomorphic to $\s(2)$, and this algebra only has inner derivations, we have in the first case that $L\otimes_k\bar k=S\otimes_k\bar k\cong\s(2,\bar k)$. 
\\
Finally, it remains to consider the case that $p=5$ and $S\otimes_k\bar k$ is isomorphic to $W(\bar k|1;\underline 1)$. We observe first, that $W(\bar k|1;\underline 1)$ has only inner derivations and is a restricted Lie algebra. Then $L\otimes_k\bar k=S\otimes_k\bar k$ and hence $L=S$. \cite[Theorem 13]{Jac} states that there is a commutative algebra ${\cal A}=k[X]/(X^5-\xi)$ for which $L\cong \Der {\cal A}$. Since $k$ is finite, it is perfect. There exists $\eta\in k$ with $\eta^5=\xi$. Set  $Y:=X-\eta$ to obtain ${\cal A}=k[Y]/(Y^5)={\cal O}(k|1;\underline 1)$. This shows that $L\cong \Der {\cal A}\cong W(k|1;\underline 1)$ and thereby completes the proof of the theorem.
\end{Proof}

\section{Lie algebras of dimension 3}

 In this section we determine all 3-dimensional  Lie algebras  over a finite field $k$ and the small dimensional irreducible modules for the simple ones.
 \\[3mm]
If $p=2$, then the perfectness of the finite field $k$ implies that $k^\ast=(k^\ast)^2:=\{\xi^2\mid \xi\in k^\ast\}$. Suppose $p>2$. 
 We observe that the homomorphism of multiplicative groups $k^\ast\rightarrow (k^\ast)^2$, $\xi\mapsto \xi^2$ has kernel $\{\pm1\}$ and therefore $|k^\ast|=2|(k^\ast)^2|$. This shows that

\begin{equation}\begin{array}{cr}
k^\ast=(k^\ast)^2\quad&\text{if }p=2,\\
\exists\delta_0 \in k^\ast\text{ so that }k^\ast=(k^\ast)^2\cup \delta_0 (k^\ast)^2\quad&\text{if } p>2.
\end{array}\end{equation} 
 
 \medskip

\begin{prop}\cite{ZP}\label{3.5}
The isomorphism classes of solvable 3-dimensional Lie algebras $R$ over the finite field $k$ are given by the following representatives.
\begin{enumerate}
\item $R$ is abelian;
\item $R=kh\oplus kx\oplus kz$, where 
$[h,x]=x,\ [z,R]=\{0\}$;
\item $R=kx\oplus ky\oplus kz$, where 
$[x,y]=z,\ [z,R]=\{0\}$;
\item $R=kd\oplus R^{(1)}$, where $R^{(1)}$ is abelian and the action of $d$ on $R^{(1)}$ is given by one of the following matrices
\begin{align*} &
 \left(\begin{array}{cc}1&0\\ 0&1
\end{array}\right), \quad \left(\begin{array}{cc}0&\xi\\1&1
\end{array}\right)\ (\xi\in k^\ast), \intertext{or}&\left(\begin{array}{cc}1&0\\ 0&-1
\end{array}\right) ,\ \left(\begin{array}{cc}0&\delta_0\\ 1&0
\end{array}\right)\text{ if }p\ne2,\quad \left(\begin{array}{cc}1&1\\ 0&1
\end{array}\right) \text{ if }p=2.\end{align*}
\end{enumerate}
\end{prop}
\begin{Proof}
(1) If $\dim R^{(1)}\le1$ only the algebras of 1., 2., 3. can occur (see for example \cite[page 34]{SF}). If $\dim R^{(1)}=2$, then $R=kd\oplus R^{(1)} $ and $ R^{(1)}$ is abelian (see for example \cite[page 34]{SF}). Thus we only have  to determine the action of $d$ on $ R^{(1)}$. Note that $d$ acts invertibly on $ R^{(1)}$, because $ R^{(1)}=[d, R^{(1)}]+[R^{(1)},R^{(1)}]=[d, R^{(1)}]$. \\
Let $\chi=T^2+\alpha T+\beta$ be the characteristic polynomial of $d$ acting on $ R^{(1)}$.  
If $\alpha=0$, we may adjust $d$ by a nonzero scalar to obtain $-\beta=1$ if $p=2$, and $-\beta\in\{1,\delta_0\}$ if $p>2$ (see equation (3.1)). If $\alpha\ne0$, we may adjust $d$ by a nonzero scalar to obtain  $\alpha=-1$. This means that we may assume
$$\chi\in\{  T^2-1,\ T^2-\delta_0,\   T^2-T+\beta,\ (\beta\in k^\ast )\}.$$
(i) Consider the case $\alpha=0$.\\
 If $p=2$, then  $\chi=T^2-1=(T-1)^2$. Then either $d={\rm Id}$ or there is $y$ such that $[d,y]-y=:x$ is linearly independent of $y$. In either case we may  choose $x,y\in R^{(1)}$ such that $[d,y]=\xi x+y$, $[d,x]=x$, $ \xi\in\{0,1\}$.\\
Suppose $p>2$. If $\chi=T^2-1$, then $\chi=(T-1)(T+1)$ and we may  choose $x,y\in R^{(1)}$ such that $[d,x]=x$, $[d,y]=-y$.\\
If $\chi=T^2-\delta_0$, then $d$ has no eigenvalue in $k$. We may  choose a basis $(x,y:=[d,x])$, and $[d,y]=(\ad d)^2(x)=\delta_0 x$ holds. \\
(ii) Consider the case $\chi=T^2-T+\beta$. 
If  there is a vector $x$ which is not an eigenvector, then  we may  choose a basis $(x,y:=[d,x])$, and $[d,y]=(\ad d)^2(x)=-\beta x+ y$ holds. Set $\xi:=-\beta$.
\\
Finally, let all nonzero vectors be eigenvectors. Then there is only one eigenvalue $r$ and $\chi=(T-r)^2$. The requested form for $\chi$ implies $p\ne2$ and $r=1/2$. Adjusting $d$ we may assume  that $d={\rm Id}$.
\\[2mm]
(2) We discuss isomorphisms between the exposed algebras. Let $\sigma:R\to R'$ be an isomorphism. For an obvious reason we only have to deal with the case that $ R^{(1)}$ and $ (R')^{(1)}$ are 2-dimensional. Let $(d,x,y)$ and $(d',x',y')$ be bases of $R$ and $R'$, respectively, for which the action of $d$ and $d'$ on $ R^{(1)}$ and $ (R')^{(1)}$ are given by matrices from the above list. \\
Note that $\sigma(d)\equiv rd' \pmod{ (R')^{(1)}}$ for some $r\in k^\ast$.    If $d$ acts as the identity, then so does $d'$ because no other matrix in the list is a multiple of the identity.
\\
Next we look at the remaining matrices. The representing matrix for $d'$ with respect to $(\sigma(x),\sigma(y))$ is $r$-times the representing matrix for $d$ with respect to $(x,y)$.
We therefore have to decide when the respective characteristic polynomials
$$ T(T-r)-r^2\xi;\quad T^2-r^2,\ T^2-r^2\delta_0\ (p\ne2);\quad T^2-r^2\ (p=2)$$
matches one of 
$$T(T-1)-\xi';\quad T^2-1,\ T^2-\delta_0\ (p\ne2);\quad T^2-1\ (p=2).$$
It is obvious that only the polynomials for the same type of matrix match and in the first case $r=1$, $\xi=\xi'$ hold.
\end{Proof}

The number of isomorphism classes of 3-dimensional solvable Lie algebras over a finite field $k$ is
$$4+|k|, \text{ if }p=2,\quad 5+|k|,\text{ if }p>2.$$
Next we determine the nonsolvable 3-dimensional Lie algebras.

\begin{theo}\label{3.1}
Let $L$ be a nonsolvable Lie algebra of dimension $\le 3$ over a finite field $k$ of characteristic $p$. \\[2mm] If $p=2$, then $L\cong W(1;\underline 2)^{(1)}$. \\[2mm]
If $p\ge 3$, then $L\cong \s(2,k)$.
\end{theo}
\begin{Proof} Recall that $L$ is 3-dimensional and simple.
\\[2mm]
(1) Consider the case $p=2$. Since $L$ is not nilpotent there is $h\in L$ which is not $\rm ad$-nilpotent. The characteristic polynomial of $\ad h$ has degree 3. Since 0 is an eigenvalue, and $tr(\ad h)=0$, this polynomial is of the form $T(T^2+\alpha)$. As $\ad h$ is not nilpotent, $\alpha\ne0$. The perfectness of $k$ shows that $\alpha=\beta^2$ is a square in $k$. Adjusting $h$ we may assume that $\ad h$ satisfies the polynomial
$T(T+1)^2$.
Choose $u,v\in L$ so that
$$L=kh\oplus ku\oplus kv,\ [h,u]=u,\ [h,v]-v\in ku.$$
Comparing eigenvalues  one gets  $[u,v]=\gamma h$, and, as $L$ is simple, $\gamma\ne0$ holds. Adjusting $u$ one may assume $\gamma=1$.
We intend to find a triple $(h',u',v')$ for which in addition $[h',v']=v'$ holds. The perfectness of $k$ allows to choose $\mu \in k$ with $\mu^2=\lambda $. Put
$$h':=h+\mu u,\ v':=\mu v+\lambda h+\mu u,\ u':=\mu ^{-1}u.$$
Then
\begin{align*}
[h',u']&=u',\\ [u',v']&=\mu^{-1}(\mu h+\lambda u)=h+\mu u=h',\\ [h',v']&= [h,v']+\mu [u,v']=\mu v+\mu \lambda u+\mu u+\mu (\mu h+\lambda u)\\&=\mu v+\mu u+\lambda h=v'.
\end{align*}
The identification 
$$v'=\partial,\ h'=x\partial,\ u'=x^{(2)}\partial$$
gives the desired isomorphism.
\\[2mm]
(2) Consider the case $p>2$.  As in the former case there is a mapping $q:L\rightarrow k$ such that the characteristic polynomial of $\ad x$ is $T(T^2+q(x))$, and there is $h\in L$ with $q(h)\ne0$. If $q(x)=0$, then $0=-q(x)$ is the only eigenvalue of $(\ad x)^2$. If $q(x)\ne0$, then one can decompose $L=kx\oplus L_1(\ad x)$ into the Fitting components with respect to $\ad x$, and  $-q(x)$ is the only eigenvalue of $(\ad x)^2$ on $L_1(\ad x)\cong L/kx$. Moreover, $T^2+q(x)$ is the characteristic polynomial of $\ad x$ on $L_1(\ad x)$ in this case.\\[2mm]
(i) We intend to show that $q$ is a quadratic form on $k^3$. Decompose $L=kh\oplus L_1$ into the Fitting components with respect to $\ad h$ (recall that $q(h)\ne0$ and $T^2+q(h)$ is the characteristic polynomial of $\ad h$ on $L_1$). If $k$ contains  a nonzero eigenvalue $\alpha$ of $\ad h$, then $-\alpha$ is also an eigenvalue (and different from $\alpha$). Therefore there is in any case  a nonzero vector $u\in L_1$ which is not an eigenvector for $\ad h$. Then $u$, $v:=[h,u]$ are linearly independent and therefore span $L_1$. Since $(\ad h)^2+q(h){\rm Id}$ vanishes on $L_1$, we have $[h,v]=-q(h)u$. From this we deduce that $[u,v]\in\ker(\ad h)=kh$. Put $[u,v]=\beta h$ with $\beta \in k$.  If $\beta =0$, then $L_1$ would be an ideal of $L$, which is not true. Therefore $\beta \ne0$. Next we compute 
$[u,[u,v]]=[u,\beta  h]=-\beta  v$, and therefore $-\beta $ is the uniquely determined nonzero eigenvalue of $(\ad u)^2$. This gives $\beta =q(u)$. The multiplication is now given by
$$[h,u]=v,\ [h,v]=-q(h)u,\ [u,v]=q(u)h.$$
For $\lambda,\kappa,\mu \in k$ we compute $q(\lambda h+\kappa u+\mu v)$ as follows
\begin{align*}
&(\ad (\lambda h+\kappa u+\mu v))^2(h)\\&=[\lambda h+\kappa u+\mu v,-\kappa v+q(h)\mu u]\\&= \lambda \kappa q(h)u+\lambda q(h)\mu v+(-\kappa^2q(u)-\mu^2q(h)q(u))h\\&= \lambda q(h)(\lambda h+\kappa u+\mu v)-(\lambda^2q(h)+\kappa^2q(u)+\mu^2q(u)q(h))h.
\end{align*}
Arguing on $L/k(\lambda h+\kappa u+\mu v)$ this gives 
$$q(\lambda h+\kappa u+\mu v)=\lambda^2q(h)+\kappa^2q(u)+\mu^2q(u)q(h).$$
Hence $q$ is a quadratic form, $q(v)=q(u)q(h)$ and $(h,u,v)$ is orthogonal with respect to the associated  bilinear form. 
\\[2mm]
(ii) Next we intend to find a nonzero element $e\in L$ with $q(e)=0$. It is a standard fact, that every anisotropic quadratic form over a finite field is at most 2-dimensional. Since the proof is very short, we give an argument to show that every quadratic form occurring in our context is isotropic. Recall the definition of $\delta_0$ from equation (3.1).
Adjusting $h,u,v$ by suitable scalars  we may substitute $q(h),q(u)$ by any elements of the same residue classes. Thus we are only interested in the following quadratic forms
$$ q=(q(h),q(u),q(u)q(h))\in\{(1,1,1),\ (1,\delta_0 ,\delta_0 ),\ (\delta_0 ,1,\delta_0 ),\ (\delta_0 ,\delta_0 ,\delta_0 ^2)\}.$$
If $-1\in(k^\ast)^2$ (so that $-1=\mu ^2$ for a suitable element $\mu \in k^\ast$), then in the respective cases we find the isotropic vector
$$e:=h+\mu u,\ u+\mu v,\ h+\mu v,\ h+\mu u.$$
If $-1\notin (k^\ast)^2$, then we may take $\delta_0 =-1$. The forms $$(1,\delta_0 ,\delta_0 ),\ (\delta_0 ,1,\delta_0 ),\ (\delta_0 ,\delta_0 ,\delta_0 ^2)$$ are isotropic, namely there is the isotropic vector 
$$e:= h+u,\ h+u,\ h+v$$
in the respective cases. It remains to consider the form  $(1,1,1)$. Suppose that for $\kappa \in \F_p$ the following implication holds
$$\kappa \in (k^\ast)^2\Rightarrow \kappa +1 \in  (k^\ast)^2.$$
Since $1\in (k^\ast)^2$ this would give the contradiction $-1\in (k^\ast)^2$. Therefore there is $\kappa \in (k^\ast)^2$ for which $\kappa +1\notin (k^\ast)^2$. Recall that $\kappa+1\ne0$ because $-1\notin (k^\ast)^2$. Then $\kappa+1\in-(k^\ast)^2$. Set $\kappa =:\lambda ^2$ and $-(\kappa +1)=:\rho ^2$. Then $1+\lambda ^2+\rho ^2=0$ and $e:=h+\lambda u+\rho v$ is isotropic.\\[2mm]
(iii) We have now constructed a nonzero element $e\in L$ satisfying $e\ne0$,  $(\ad e)^3=0$.  Clearly, $\ker(\ad e)$ is not 3-dimensional (otherwise $e\in C(L)=\{0\}$). Suppose it is 2-dimensional. Then $\ker(\ad e)=ke\oplus kx$ is abelian, and $[e,L]\subset \ker(\ad e)$ (as $\ad e$ is nilpotent). Choose $y\not\in \ker(\ad e)$. We obtain
$$[e,[x,y]]=[x,[e,y]]\in [x,\ker(\ad e)]=\{0\}.$$
But then $[L,\ker(\ad e)]\subset \ker(\ad e)$, and $\ker(\ad e)$ would be an ideal. This contradiction shows $\dim\, \ker(\ad e)=1$, whence   $\ker(\ad e)=ke$. Next suppose that $(\ad e)^2=0$. Then $[e,L]\subset ke$, which implies that $\dim\, L=\dim\, [e,L]+\dim\,\ker(\ad e)\le 2$, a contradiction. Hence $(\ad e)^2(L)=ke$.
Choose $h\in [e,L]$ so that $[h,e]=2e$. Next find $f\in L$ with $h=[e,f]$. We now observe that this can only be if $0,\pm2$ are eigenvalues of $\ad h$. But then all the eigenspaces are 1-dimensional, and we may take $f$ as an eigenvector for the eigenvalue $-2$. Consequently, $L\cong \s(2,k)$.
\end{Proof}

Next we are interested in at most 3-dimensional modules of these simple 3-dimensional Lie algebras. Recall that for $p=2$ the algebra $W(1;\underline 2)^{(1)}$ is simple with basis $(\partial, x\partial, x^{(2)}\partial)$, and ${\cal O}(1;\underline 2)/k$ is a 3-dimensional module for this algebra with basis $(x+k, x^{(2)}+k, x^{(3)}+k)$.

\begin{prop}\label{3.2}
Let $p=2$.
\begin{enumerate}
\item $\Der W(1;\underline 2)^{(1)}=k\partial^2\oplus W(1;\underline 2)$.
\item Let $\alpha,\beta,\gamma,\delta\in k$ satisfy $\alpha\delta+\beta\gamma=1$. The mapping
\begin{align*}
&\partial\mapsto \alpha\partial+\beta x^{(2)}\partial,\quad x\partial\mapsto x\partial,\quad  x^{(2)}\partial\mapsto \gamma\partial+\delta x^{(2)}\partial,\\ &
\partial^2\mapsto \alpha^2\partial^2+\alpha\beta x\partial+\beta^2 x^{(3)}\partial,\quad  x^{(3)}\partial \mapsto \gamma^2\partial^2+\gamma\delta x\partial+\delta^2 x^{(3)}\partial
\end{align*}
defines an automorphism of $\Der W(1;\underline 2)^{(1)}$.
\item Every faithful $W(1;\underline 2)^{(1)}$-module has dimension bigger than 2.
\item The only   faithful $W(1;\underline 2)^{(1)}$-module of dimension 3 is ${\cal O}(1;\underline 2)/k$.
\end{enumerate}
\end{prop}
\begin{Proof}
(1) Let $d$ denote a derivation of $W(1;\underline 2)^{(1)}$. Adjusting $d$ by adding  suitable multiples of $x\partial$, $x^{(2)}\partial$ and $x^{(3)}\partial$ we may assume that $d(\partial)=0$. Adding  suitable multiples of $\partial$ and $\partial^2$ we may in addition assume that $d(x^{(2)}\partial)\in k x^{(2)}\partial$. Set $d(x^{(2)}\partial)=\alpha x^{(2)}\partial$ for some $\alpha\in k$. Then  $0=d(\partial)=d([\partial,[\partial,x^{(2)}\partial]])=[\partial,[\partial,d(x^{(2)}\partial)]]=\alpha \partial$. This gives $\alpha=0$  and $d=0$.
\\[2mm]
(2) easy computation.
\\[2mm]
(3) Since $p=2$, one has that $\frak{gl}(2)$ is solvable. But $W(1;\underline 2)^{(1)}$ is simple and therefore does not fit into $\frak{gl}(2)$.
\\[2mm]
(4) Set $L:=W(1;\underline 2)^{(1)}$, let $V$ be a faithful 3-dimensional $L$-module and $\rho: L\to \frak{gl}(V)$ the representation. By (3) $V$ is an irreducible module. Recall that ${\cal C}(V,L)$ is a field extension of $k$ and $V$ is a ${\cal C}(V,L)$-vector space. For dimension reasons one obtains hat either $V={\cal C}(V,L)v$ for some $v\in V$ or ${\cal C}(V,L)=k$ is true. In the first case there is a  mapping $\Phi:L\to {\cal C}(V,L)$ such that $\rho(x)(\alpha v)=\alpha \rho(x)(v)=\alpha \Phi(x)v=\Phi(x)\alpha v$ for all $\alpha\in {\cal C}(V,L)$. But then $L$ would be abelian, which is not true. \\
Therefore we have ${\cal C}(V,L)=k$. Note that $\rho(\partial)^4$, $\rho(x\partial)^2+\rho(x\partial)$ and $\rho(x^{(2)}\partial)^4$ are contained in ${\cal C}(V,L)$. Therefore there are $\alpha,\beta,\gamma\in k$ for which
$$\rho(\partial)^4=\alpha{\rm Id},\ \rho(x\partial)^2+\rho(x\partial)=\beta{\rm Id},\ \rho(x^{(2)}\partial)^4=\gamma{\rm Id}.$$
On the other hand, one has $tr(\rho(x))=0$ for every $x\in W(1;\underline 2)^{(1)}$ (since this algebra is simple), and as $tr(\rho(x)^{2^i})=(tr(\rho(x)))^{2^i}$ for every $x$ whereas $tr({\rm Id}_V)=3\ne0$, this gives $\alpha=\beta=\gamma=0$. Then $\rho(\partial)^3=0$ and there is a vector $v\in V$ for which
$$\rho(x^{(2)}\partial)(v)=0,\ \rho(x\partial)(v)=\delta v\text{ for some }\delta\in \F_2.$$
The irreducibility of $V$ yields 
$$V=\oplus_{i=0}^2 k\rho(\partial)^i(v).$$
If $\delta=0$, then an easy computation shows that $\oplus_{i=1}^2 k\rho(\partial)^i(v)$ is a proper submodule. Therefore this case is impossible. Then $\delta=1$, and this shows that there is only one 3-dimensional irreducible module. This module then has to be as claimed.
\end{Proof}

 For $p\ge3$ the algebra $\s(2,k)$ plays a distinguished role in our context.  Let $(e,h,f)$ be a basis for $\s(2,k)$ which we call an {\em $\s(2)$-triple} if
 \begin{equation}[h,e]=2e,\ [h,f]=-2f,\ [e,f]=h.\end{equation}
 Note that $\s(2,k)$ is a restricted Lie algebra. If $V$ is an irreducible restricted module, then $e$ acts nilpotently on $V$. The standard procedure shows  there is a vector $v_0\in V$ for which 
 $$\rho(e)(v_0)=0,\ \rho(h)(v_0)=\alpha v_0,\ \alpha=\dim\,V-1,\ \text{and }V=\oplus_{i=0}^\alpha k\rho(f)^i(v_0) .$$
Also, $\dim V\le p$ holds.  These modules  are denoted by $V(\alpha)$.  In general, $\rho(x)^p-\rho(x^{[p]})$ is, for all $x\in \s(2,k)$,  contained in ${\cal C}(V,L)$.  It is well known that, denoting the algebraic closure of $k$ by $\bar k$, there is  a linear form $\chi:\s(2,k)\to \bar k$ such that 
  $$\rho(x)^p-\rho(x^{[p]})=\chi(x)^p {\rm Id}\quad\forall\ x\in\s(2,k).$$
The linear form $\chi$ is called the character of the representation.  We call a $k$-linear form
$\chi:\s(2,k)\to  k$
an {\em irreducible character of dimension $d$}, if there is an irreducible $d$-dimensional representation $\rho:\s(2,k)\to \frak{gl}(V)$ such that $\rho(x)^p=\rho(x^{[p]})+\chi(x)^p {\rm Id}$ for all $ x\in\s(2,k)$.

\begin{prop}\label{3.3} Let $k$ be any field of characteristic $3$ and $(e,h,f)$ an $\s(2)$-triple.
\begin{enumerate}
\item An element $x=\alpha e+\beta h+\gamma f$ is $\ad$-nilpotent if 
$\alpha\gamma+\beta^2=0$
and toral if 
$\alpha\gamma+\beta^2=1$.
\item Let $\sigma_{\alpha,\beta}\in {\rm GL}(\s(2,k))$, $\alpha\in k^\ast$, $\beta\in k$ be defined by
$$\sigma_{\alpha,\beta}(e):=\alpha f,\ \sigma_{\alpha,\beta}(h):=-h+\alpha\beta f,\ \sigma_{\alpha,\beta}(f):=\alpha^{-1}e-\beta h-\alpha\beta^2 f.$$
\begin{enumerate}
\item Every $\sigma_{\alpha,\beta}$ is an automorphism of $\s(2,k)$.
\item $\sigma_{\alpha,0}^2={\rm Id},\quad \sigma_{1,0}\circ \exp(\ad \beta e)=\sigma_{1,\beta},\quad \sigma_{\alpha,0}\circ \sigma_{\beta,0}\circ\sigma_{\gamma,\delta}=\sigma_{\alpha\beta^{-1}\gamma,\delta}$.
\item Every automorphism of $\s(2,k)$ is of the form
$$\sigma_{\gamma,\delta}\quad \text{or}\quad \sigma_{1,\beta}\circ\sigma_{\gamma,\delta}\quad \text{for some }\gamma\in k^\ast,\ \beta,\delta\in k.$$
\end{enumerate}
\end{enumerate}
\end{prop}
\begin{Proof}
(1) Note that (see \cite[page 64]{SF})
\begin{align*}
x^{[3]}&=(\alpha e+\beta h+\gamma f)^{[3]}\\&=\alpha^3e^{[3]}+[\alpha e,[\alpha e,\beta h+\gamma f]]+[\beta h+\gamma f,[\beta h+\gamma f,\alpha e]]+(\beta h+\gamma f)^{[3]}\\
&=[\alpha e,\alpha\beta e+\alpha\gamma h]+[\beta h+\gamma f,-\alpha\beta e-\alpha\gamma h]\\&\qquad +\beta^3h^{[3]}+[\beta h,[\beta h,\gamma f]]+[\gamma f,[\gamma f,\beta h]]+\gamma^3f^{[3]}\\&=\alpha^2\gamma e+\alpha\beta^2 e+\alpha\beta\gamma h+\alpha\gamma^2 f+\beta^3 h+\beta^2\gamma f
\\&=(\alpha\gamma+\beta^2)\alpha e+(\alpha\gamma+\beta^2)\beta h+(\alpha\gamma+\beta^2)\gamma f\\
&=(\alpha\gamma+\beta^2)x.
\end{align*}
 This proves the claim.
\\[2mm]
(2) (i) To show that $\sigma_{\alpha,\beta}$ is an automorphism we compute
\begin{align*}
[\sigma_{\alpha,\beta}(e),\sigma_{\alpha,\beta}(h)]&=[\alpha f,-h+\alpha\beta f]=\alpha f=\sigma_{\alpha,\beta}(e),\\
[\sigma_{\alpha,\beta}(e),\sigma_{\alpha,\beta}(f)]&=[\alpha f,\alpha^{-1}e-\beta h-\alpha\beta^2f]=-h+\alpha\beta f=\sigma_{\alpha,\beta}(h),\\
[\sigma_{\alpha,\beta}(h),\sigma_{\alpha,\beta}(f)]&=[-h+\alpha\beta f,\alpha^{-1}e-\beta h-\alpha\beta^2f]\\&=\alpha^{-1}e+\alpha\beta^2f-\beta h+\alpha\beta^2f=\sigma_{\alpha,\beta}(f).
\end{align*}
(ii)
The equation $\sigma_{\alpha,0}^2={\rm Id}$ is obviously true.
\\[2mm]
Next, $\sigma_{1,0}\circ \exp(\ad \beta e)$ maps 
\begin{align*}
e&&\text{onto}&&&\sigma_{1,0}(e)=\sigma_{1,\beta}(e),\\
h&&\text{onto}&&&\sigma_{1,0}(h+\beta e)=-h+\beta f=\sigma_{1,\beta}(h),\\
f&&\text{onto}&&&\sigma_{1,0}(f+\beta h-\beta^2 e)=e-\beta h-\beta^2 f=\sigma_{1,\beta}(f);
\end{align*}
next, $\sigma_{\alpha,0}\circ\sigma_{\beta,0}$ maps
$$e\mapsto \alpha^{-1}\beta,\ h\mapsto h,\ f\mapsto \alpha\beta^{-1}f$$
and therefore 
$\sigma_{\alpha,0}\circ\sigma_{\beta,0}\circ\sigma_{\gamma,\delta}$ maps
\begin{align*}
e&&\text{onto}&&&\alpha\beta^{-1}\gamma e=\sigma_{\alpha\beta^{-1}\gamma,\delta}(e),\\
h&&\text{onto}&&&-h+\alpha\beta^{-1}\gamma\delta f =\sigma_{\alpha\beta^{-1}\gamma,\delta}(h),\\
f&&\text{onto}&&& \alpha^{-1}\beta\gamma^{-1}e-\delta h-\alpha\beta^{-1}\gamma\delta^2 f=\sigma_{\alpha\beta^{-1}\gamma,\delta}(f).
\end{align*}
(iii) Let $\sigma$ be an automorphism of $\s(2,k)$. Then $\sigma(e)$ is $\ad$-nilpotent.\\
 Consider first the case that $\sigma(e)=e$:
write $\sigma(h)=\alpha e+\beta h+\gamma f$, $\sigma(f)=\alpha' e+\beta' h+\gamma' f$. Then
\begin{align*}
e=\sigma(e)&=\sigma([e,h])=[e,\sigma(h)]=\beta e+\gamma h,\\
\alpha e+\beta h+\gamma f=\sigma(h)&=\sigma([e,f])=[e,\sigma(f)]=\beta'e+\gamma'h.
\end{align*}
This gives $\beta=1$, $\gamma=0$, $\beta'=\alpha$, $\gamma'=\beta$. Furthermore, 
$$
\alpha' e+\beta' h+\gamma' f=\sigma(f)=\sigma([h,f])=[\sigma(h),\sigma(f)]=\alpha^2e+\alpha h-\alpha' e+f,$$
and this gives $ 2\alpha'=\alpha^2$. As a result, $\sigma=\exp(\ad(\alpha e))=\sigma_{1,0}\circ \sigma_{1,\alpha}$ by (ii).

Now consider the general case: write $\sigma(e)=\alpha e+\beta h+\gamma f$. If $\alpha=0$, then (1) implies $\beta=0$. In this case $\sigma(e)=\gamma f=\sigma_{\gamma,0}(e)$, and according to the previous case we have for some $\delta\in k$ (applying (ii)) $\sigma_{\gamma,0}^{-1}\circ\sigma=\sigma_{1,0}\circ\sigma_{1,\delta}$. Then
$$\sigma=\sigma_{\gamma,0}\circ\sigma_{1,0}\circ\sigma_{1,\delta}=\sigma_{\gamma,\delta}.$$
If $\alpha\ne0$, then (1) implies $\gamma=-\alpha^{-1}\beta^2$, whence 
\begin{align*}\sigma(e)&=\alpha e+[-\alpha^{-1}\beta f,\alpha e]-[-\alpha^{-1}\beta f,[-\alpha^{-1}\beta f,\alpha e]]\\&=\exp(\ad(-\alpha^{-1}\beta f))(\alpha e)=(\exp(\ad(-\alpha^{-1}\beta f))\circ\sigma_{1,0})(\alpha f)\\&=(\sigma_{1,0}\circ\exp(\ad(-\alpha^{-1}\beta e)))(\alpha f)=\sigma_{1,-\alpha^{-1}\beta}(\alpha f).
\end{align*}
The previous case gives for some $\delta\in k$
$$\sigma=\sigma_{1,-\alpha^{-1}\beta}\circ\sigma_{\alpha,\delta}.$$
\end{Proof}

\begin{prop}\label{3.4} Let $k$ be a finite field of characteristic $3$ and $(e,h,f)$ an $\s(2)$-triple. For any linear mapping $\chi:\s(2,k)\to k$ write $\chi=(\chi(e),\chi(h),\chi(f))$. The ${\rm Aut}\,\s(2,k)$-orbits of the set of nonzero irreducible 3-dimensional characters are in 1-1-corrrespondence with
$$(1,0,\xi)\text{ where }\xi\in k \text{ allows a solution of }T^3+T^2=\xi\text{ in }k.$$
\end{prop}
\begin{Proof}
(1) Note that ${\rm Aut}\,\s(2,k)$ acts on the space of all linear forms by $(\sigma\cdot\chi)(x)=\chi(\sigma^{-1}(x))$. In order to determine orbits we have to compute (cf. Proposition \ref{3.3}) the following
\begin{align*}
\sigma_{\alpha,\beta}(e)&=\alpha f,\\
\sigma_{\alpha,\beta}(h)&=-h+\alpha\beta f,\\
\sigma_{\alpha,\beta}(f)&=\alpha^{-1}e-\beta h-\alpha\beta^2f;\\[2mm]
\sigma_{1,\beta}\circ\sigma_{\gamma,\delta}(e)&=\sigma_{1,\beta}(\gamma f)=\gamma(e-\beta h-\beta^2f),\\
\sigma_{1,\beta}\circ\sigma_{\gamma,\delta}(h)&=\sigma_{1,\beta}(-h+\gamma\delta f)=-(-h+\beta f)+\gamma\delta(e-\beta h-\beta^2f),\\&=\gamma\delta e+(1-\beta\gamma\delta)h-(\beta+\beta^2\gamma\delta)f,\\
\sigma_{1,\beta}\circ\sigma_{\gamma,\delta}(f)&=\sigma_{1,\beta}(\gamma^{-1}e-\delta h-\gamma\delta^2f)\\&= \gamma^{-1}f-\delta(-h+\beta f)-\gamma\delta^2(e-\beta h-\beta^2f)\\&=-\gamma\delta^2e+(\delta+\beta\gamma\delta^2)h+(\gamma^{-1}-\beta\delta+\beta^2\gamma\delta^2)f.
\end{align*}
Therefore the orbit of a nonzero character $\chi=(r,s,t)$ consists exactly of all characters
\begin{align*}
&(\alpha t,-s+\alpha\beta t,\alpha^{-1}r-\beta s-\alpha\beta^2t),\hspace{4cm} \alpha\in k^\ast,\ \beta\in k,\\
&(\gamma (r-\beta s-\beta^2 t), \gamma\delta r+(1-\beta\gamma\delta)s-(\beta+\beta^2\gamma\delta)t,\\
&\hspace{2cm}-\gamma\delta^2r+(\delta+\beta\gamma\delta^2)s+(\gamma^{-1}-\beta\delta+\beta^2\gamma\delta^2)t),\  \gamma\in k^\ast,\ \beta,\delta\in k.
\end{align*}
Since one of $r,s,t$ is nonzero, there is a choice of $\alpha,\beta$ for which  $\alpha^{-1}r-\beta s-\alpha\beta^2t\ne0$. Therefore  the ${\rm Aut}\,\s(2,k)$-orbit of a nonzero character $\chi$ contains a character $\chi'=(r',s',t')$ with $t'\ne0$. Then we choose $\alpha:=t'^{-1}$, $\beta:=s'$ and obtain a character $(1,0,\xi)$ in this orbit.\\[2mm]
(2) Next we determine when characters $(1,0,\xi)$ and $(1,0,\xi')$ are in the same orbit. \\
If $(1,0,\xi')=\sigma_{\alpha,\beta}^{-1}\cdot (1,0,\xi)$, then (1) implies
$$\alpha\xi=1,\ \alpha\beta\xi=0,\ \alpha^{-1}-\alpha\beta^2\xi=\xi'.$$
This gives $\beta=0$ and $\xi'=\alpha^{-1}=\xi$.
\\
 If $(1,0,\xi')=(\sigma_{1,\beta}\circ\sigma_{\gamma,\delta})^{-1}\cdot (1,0,\xi)$, then (1) gives
 $$\gamma(1-\beta^2\xi)=1,\ \gamma\delta-(\beta+\beta^2\gamma\delta)\xi=0,\ -\gamma\delta^2+(\gamma^{-1}-\beta\delta+\beta^2\gamma\delta^2)\xi=\xi'.$$
 These equations imply
 \begin{align*}
&\gamma(1-\beta^2\xi)=1,\qquad \beta\xi=\delta,\\
&\xi'=-\gamma\delta^2(1-\beta^2\xi)+(1-\beta^2\xi)\xi-\delta^2\\&\ =-\delta^2+\xi-\beta^2\xi^2-\delta^2=\xi
\end{align*}
(3) We have to determine those $\xi$ for which there is a 3-dimensional irreducible $\s(2,k)$-module $V$ having character $(1,0,\xi)$. Let $V$ be such a module with representation $\rho$. Then
$$ \rho(e)^3={\rm Id},\ \rho(h)^3=\rho(h),\ \rho(f)^3=\xi^3{\rm Id}.$$
 We observe that $\rho(h)$ is a semisimple endomorphism and has all eigenvalues contained  in $\F_3$. Choose an eigenvector $v_0$ with $\rho(h)$-eigenvalue $\mu\in \F_3$. Since $\rho(e)$ is invertible, one has that $v_0$, $\rho(e)v_0$, $\rho(e)^2v_0$ are nonzero eigenvectors for $\rho(h)$ with respective eigenvalues $\mu$, $2+\mu$, $1+\mu$. Thus these vectors are linearly independent and span $V$. Substituting $v_0$ by one of the others we may assume $\mu=0$. Considering eigenvalues one obtains that $\rho(f)v_0=\lambda \rho(e)^2v_0$ for some $\lambda\in k$. One computes
 \begin{align*}
\rho(f)\rho(e)v_0&=-\rho(h)v_0+\lambda\rho(e)^3v_0=\lambda v_0,\\
\rho(f)\rho(e)^2v_0&=-\rho(h)\rho(e)v_0+\rho(e)\rho(f)\rho(e)v_0=(-2+\lambda)\rho(e)v_0.\end{align*}
Then $\xi^3{\rm Id}=\rho(f)^3=\lambda^2(-2+\lambda){\rm Id}$. Since $k$ is perfect, one finds $\kappa\in k$ satisfying $\kappa^3=\lambda$. Then $\kappa$ solves the equation $T^3+T^2=\xi$.
\\
Conversely, if the equation in question is solvable with solution $\kappa$, then the above exposed action defines a module action, and the module clearly is irreducible.
\end{Proof}

\section{Nonsolvable Lie algebras of dimension 4 and 5}

In this chapter we determine the algebras mentioned in the title.

\begin{theo}\label{4.1}Let $L$ be a nonsolvable Lie algebra of dimension $4$ over a finite field $k$ of characteristic $p$.  \\ If $p=2$, then $$L\in\{ W(1;\underline 2),\ W(1;\underline 2)^{(1)}\oplus C(L)\}.$$ If $p\ge 3$, then $$L\cong \frak{gl}(2,k).$$ 
\end{theo}
\begin{Proof}Note that $\dim \rad(L)\le 1$.\\[2mm]
(1) Consider the case $p=2$. 
\\[2mm]
If $\rad(L)=\{0\}$, then $L$ is semisimple. Theorem \ref{2.4} in combination with Theorem \ref{3.1} shows that $L\subset \Der W(1;\underline 2)^{(1)}$, whence $L=kd\oplus W(1;\underline 2)^{(1)}$ for some $d\in k\partial^2+kx^{(3)}\partial$ (Proposition \ref{3.2}). There is nothing to prove if $d\in kx^{(3)}\partial$. Otherwise put $d= \alpha'\partial^2+\beta'x^{(3)}\partial$, choose $\alpha,\beta,\gamma,\delta$ according to $$\gamma^2=\alpha',\ \delta^2=\beta', \alpha=0,\ \beta\gamma=1.$$
The automorphism $\sigma$ mentioned in Proposition \ref{3.2}(2) coming with this choice maps $x^{(3)}\partial-\gamma\delta x\partial$ onto $d$. Then $\sigma^{-1}(d)=x^{(3)}\partial-\gamma\delta x\partial$ and $\sigma^{-1}(L)=W(1;\underline 2)$.
\\[2mm]
If $\dim \rad(L)=1$, then $L/\rad(L)\cong W(1;\underline 2)^{(1)}$ (Theorem \ref{3.1}). The perfectness of this algebra shows that $\rad(L)=C(L)$. Choose an inverse image $h'$ of $x\partial$ and inverse images $e,f$ of $\partial$ and $x^{(2)}\partial$ which are $(\ad h')$-eigenvectors, respectively, and set $h:=[e,f]$. Note that $\ad h=\ad h'$. It is easily seen that $ke+kf+kh$ is a 3-dimensional subalgebra, hence an ideal of $L$. It is clear that $ke+kf+kh\cong W(1;\underline 2)^{(1)}$.\\[2mm]
(2) Consider the case $p\ge3$. Theorem \ref{2.4} shows that there are no semisimple 4-dimensional Lie algebras, and therefore $L$ is an extension of $\s(2,k)$ by a 1-dimensional center. This extension splits (by the same argument used in (1)).
\end{Proof}

\begin{theo}\label{4.2}Let $L$ be a nonsolvable Lie algebra of dimension $5$ over a finite field $k$ of characteristic $p=2$.  Then $L$ is one of the following. 
\begin{enumerate}
\item $\Der W(1;\underline 2)^{(1)}$;
\item $W(1;\underline 2)\rtimes \rad(L)$, where  $[W(1;\underline 2)^{(1)},\rad(L)]=\{0\}$ and $$\rad(L)=ku,\quad [x^{(3)}\partial,u]=\delta u,\qquad \delta\in\{0,1\};$$
\item $W(1;\underline 2)^{(1)}\oplus \rad(L)$ is the direct sum of ideals, and $$  \rad(L)=kh\oplus ku,\quad  [h,u]=\delta u,\qquad \delta\in\{0,1\}.$$
\end{enumerate}
The exposed Lie algebras are mutually nonisomorphic.
\end{theo}
\begin{Proof} (1) If $\rad(L)=\{0\}$, then  Theorem \ref{2.4} in combination with Theorem \ref{3.1} shows that $L\subset \Der W(1;\underline 2)^{(1)}$, while Proposition \ref{3.2} proves that $\Der W(1;\underline 2)^{(1)}$ is 5-dimensional. Then
 $L=\Der W(1;\underline 2)^{(1)}$.
\\[2mm]
(2) Suppose $\dim \rad(L)=1$. Set $\rad(L)=ku$. Theorem \ref{4.1} shows that $L/\rad(L)\cong W(1;\underline 2)$. Let $Q$ denote the inverse image  in $L$ of the subalgebra  $W(1;\underline 2)^{(1)}$. This is a 4-dimensional ideal containing the radical of $L$, whence by Theorem \ref{4.1}
$$Q\cong W(1;\underline 2)^{(1)}\oplus C(Q),\quad C(Q)=\rad(L).$$
Write $L=Q\oplus kd$, observe that $Q^{(1)}\cong W(1;\underline 2)^{(1)}$. Then $Q^{(1)}+kd$ is a 4-dimensional subalgebra of $L$. If it would have a nontrivial center, then we would have $\dim\rad(L)=2$. As this is not the case, Theorem \ref{4.1} shows that $Q^{(1)}+kd\cong W(1;\underline 2)$.
The case that $x^{(3)}\partial $ acts trivially on $\rad(L)$ is listed in the theorem. Otherwise $[x^{(3)}\partial,u]=\lambda u\ne0$. Set in Proposition \ref{3.2} $\beta=\gamma=0$, $\delta^2=\lambda^{-1}$, $\alpha\delta=1$. The automorphism coming with this choice maps $x^{(3)}\partial$ onto $\lambda^{-1}x^{(3)}\partial$. So we may assume $\lambda=1$.
\\[2mm]
(3) Suppose $\dim \rad(L)=2$. Theorem \ref{3.1} shows that $L/\rad(L)\cong W(1;\underline 2)^{(1)}$. If $\rad(L)$ is abelian, then Proposition \ref{3.2}(3) proves that $\rad(L)$ is a trivial $L/\rad(L)$-module. If $\rad(L)$ is not abelian, then $\rad(L)=kh\oplus ku$ with $[h,u]=u$. 
In both cases $L$ has a 1-dimensional ideal $ku$. Applying Theorem \ref{4.1} twice we see that $L/ku\cong W(1;\underline 2)^{(1)}\oplus kz$ splits and then that the inverse image of $W(1;\underline 2)^{(1)}$ in $L$ also splits. Hence $L$ has an ideal $P$ isomorphic to $W(1;\underline 2)^{(1)}$.
\\[2mm]
(4) In all listed cases $L$ has a unique perfect ideal $L^{(2)}\cong W(1;\underline 2)^{(1)}$, and the algebras in question are distinguished by $\dim\rad(L)$ and $L/L^{(2)}$.
\end{Proof}

\begin{theo}\label{4.3}Let $L$ be a nonsolvable Lie algebra of dimension $5$ over a finite field $k$ of characteristic $p\ge3$. The following occurs. 
\begin{enumerate}
\item $p=5$ and $ L\cong W(1;\underline 1)$;
\item $L=P\rtimes \rad(L)$, where $P\cong\s(2,k)$,  and 
\begin{enumerate}
\item $\rad(L)=C(L)$, or
\item $\rad(L)=kh\oplus ku$ with $[h,u]=u$, $[P,\rad(L)]=\{0\}$, or
\item $\rad(L)$ is abelian, $\rad(L)\cong V(1)$ is the irreducible 2-dimensional $P$-module;
\end{enumerate}
\item $p=3$, and $L$ is a nonsplit extension $0\to V(1)\to L\to \s(2,k)\to 0$; more exactly, $L$ has a basis $(e,h,f,v_0,v_1)$ with multiplication
\begin{align*}
&[h,e]=-e+v_1,&&[h,v_0]=v_0,&& [e,v_0]=0,&&[f,v_0]=v_1,\\
&[h,f]=f,&&[h,v_1]=-v_1,&&[e,v_1]=v_0,&& [f,v_1]=0,\\
&[e,f]=h,&&[v_0,v_1]=0.
\end{align*}
\end{enumerate}
The exposed algebras are mutually nonisomorphic.
\end{theo}
\begin{Proof} (1) Theorem \ref{2.4} shows that there is exactly one  5-dimensional semisimple Lie algebra, namely $W(1;\underline 1)$ in case $p=5$.
\\[2mm]
 (2) Consider the case that $\rad(L)\ne\{0\}$. There is no semisimple Lie algebra of dimension 4 (Theorem \ref{4.1}). Therefore  $\dim\rad(L)=2$ holds. Suppose $\rad(L)$ contains a 1-dimensional ideal $ku$ of $L$. Applying Theorem \ref{4.1} twice we see that $L/ku\cong \frak{gl}(2,k)$ splits and then that the inverse image of $\s(2,k)$ in $L$ also splits. Hence $L$ has an ideal isomorphic to $\s(2,k)$, which annihilates $\rad(L)$. Being a 2-dimensional algebra $\rad(L)$ is of the required form.
  \\[2mm]
 (3) Next suppose that $\rad(L)$ is $L$-simple. Then $\rad(L)$ is abelian, and $L/\rad(L)$ acts on $\rad(L)$ as an $\s(2,k)$ (Theorem \ref{3.1}). There is only one isomorphism class of irreducible $\s(2,k)$-modules of dimension 2 in characteristic $p>2$, given by $V(1)=kv_0\oplus kv_1$ and 
 $$ [h,v_0]=v_0,\ [h,v_1]=-v_1,\ [e,v_0]=0,\ [f,v_0]=v_1,\ [e,v_1]=v_0,\ [f,v_1]=0.$$
If the extension $0\to V(1)\to L\to \s(2,k)\to 0$ splits,  then we are in case 2.(c) of the theorem. So assume that the extension does not split. Choose an $\s(2)$-triple $(\bar e,\bar h,\bar f)$ in $L/\rad(L)$ and let $h$ be an inverse image of $\bar h$ in $L$. If $p>3$, then $L$ decomposes into 1-dimensional $h$-root spaces
$$L=L_{-2}\oplus L_{-1}\oplus L_0\oplus L_1\oplus L_2,\ \rad(L)=L_{-1}\oplus L_1,$$
and $L_{-2}\oplus L_0\oplus L_2$ is a subalgebra isomorphic to $\s(2,k)$. But then the extension splits, which is not true in the present case. Therefore $p=3$, and $L$ decomposes
$$L=L_{-1}\oplus L_0\oplus L_1,\quad L_0=kh,\ v_1\in L_{-1},\ v_0\in L_1.$$
Choose inverse images $e\in L_{-1}$ and $f\in L_1$ for $\bar e$ and $\bar f$, respectively. Then there are $\alpha,\beta\in k$ such that 
$$[h,e]=-e+\alpha v_1,\quad [h,f]=f+\beta v_0.$$
The nonsplitting of the sequence means $\alpha\ne0$ or $\beta\ne0$. Intertwining $e,f$ and $v_0,v_1$ if necessary gives $\alpha\ne0$ in any case. If $\beta=0$, then we set $v_0':=\alpha v_0$, $v_1':=\alpha v_1$.
\\
Otherwise we choose $\lambda\in k$ with $\lambda^3=\alpha\beta^{-1}$ and set
\begin{align*}&e':=\lambda f+\lambda \beta v_0-\lambda^2\beta v_1,\quad h':=-h+\lambda f,\quad f':=\lambda^{-1}e-h-\lambda f+\lambda\beta v_0,\\& v_1':=\lambda^{-1} v_0 +v_1,\quad v_0':=v_1.
\end{align*}
Then 
\begin{align*}
[e',f']&=[f,e]-\lambda [f,h]+\lambda^2\beta [f,v_0]+\beta [v_0,e]-\lambda\beta[v_0,h]-\lambda^2\beta[v_0,f]\\
&\quad -\lambda\beta[v_1,e]+\lambda^2\beta[v_1,h]+\lambda^3\beta[v_1,f]
\\&=-h+\lambda(f+\beta v_0)+\lambda^2\beta v_1+\lambda\beta v_0+\lambda^2\beta v_1+\lambda\beta v_0+\lambda^2\beta v_1\\&=-h+\lambda f=h',\\[2mm]
[h',e']&=-\lambda [h,f]-\lambda\beta[h,v_0]+\lambda^2\beta[h,v_1]+\lambda^2\beta[f,v_0]\\&=-\lambda(f+\beta v_0)-\lambda\beta v_0-\lambda^2\beta v_1+\lambda^2\beta v_1\\&=(-\lambda f-\lambda\beta v_0+\lambda^2\beta v_1)-(\lambda\beta v_0+\lambda^2\beta v_1)\\&=-e'-\lambda^2\beta v_1',\\[2mm]
[h',f']&= -\lambda^{-1}[h,e]+\lambda[h,f]-\lambda\beta[h,v_0]+[f,e]-\lambda[f,h]+\lambda^2\beta[f,v_0]\\&=-\lambda^{-1}(-e+\alpha v_1)+\lambda(f+\beta v_0)-\lambda\beta v_0-h +\lambda(f+\beta v_0)+\lambda^2\beta v_1\\&= \lambda^{-1} e-\lambda^{-1}(\lambda^3\beta) v_1+\lambda f+\lambda\beta v_0-\lambda\beta v_0-h+\lambda f+\lambda\beta v_0+\lambda^2\beta v_1\\&= \lambda^{-1}e+2\lambda f-h+\lambda\beta v_0=f'.
\end{align*}
This shows that $(\bar{e'},\bar{h'},\bar{f'})$ is an $\s(2)$-triple, and with this choice we are in the former case.\\[2mm]
(4) The properties of the radical distinguishes all algebras mentioned in cases 1 and  2. It remains to show that the algebra mentioned in case 3 is in fact nonsplit. So assume on the contrary that it is isomorphic to an algebra $P\rtimes V(1)$ where $P\cong \s(2,k)$. Let $(e',h',f')$ denote the projections of $(e,h,f)$ into $P$ with respect to this decomposition. The multiplication of $L/\rad (L)\cong P$ shows that  $(e',h',f')$ is an $\s(2)$-triple in $P$. Write 
\begin{align*}
e&=e'+\alpha_0 v_0+\alpha_1 v_1,\ h=h'+\beta_0 v_0+\beta_1 v_1,\ f=f'+\gamma_0 v_0+\gamma_1v_1.
\end{align*}
Since $\ad(\beta_0v_0-\beta_1v_1)^{(p+1)/2}=0$, we have that $\exp(\ad(\beta_0v_0-\beta_1v_1))$ is an automorphism of $L$. One has $\exp(\ad(\beta_0v_0-\beta_1v_1))(h)=h-\beta_0v_0-\beta_1v_1=h'$. We therefore may assume that $h=h'$. Then
$$-e+v_1=[h,e]=[h',e']+\alpha_0v_0-\alpha_1v_1=-e'+\alpha_0v_0-\alpha_1v_1=-e+2\alpha_0v_0,$$
a contradiction.
\end{Proof}

 As a result, there are exactly the following numbers of isomorphism classes of nonsolvable Lie algebras over a finite field 

\vspace{1cm}
\begin{tabular}{l|c|c|cc}
&$\dim L=3$&$\dim L=4$&$\dim L=5$\\ \hline
$p=2$&1&2&5\\$p=3,5$&1&1&4\\$p>5$&1&1&3&.
\end{tabular}

\vspace{2cm}

In order to determine the 6-dimensional Lie algebras in the next section we have to derive some subsidiary results.

\begin{prop}\label{4.4} Let $p=3$, and $L_1$ be the nonsplit extension $0\to V(1)\to L\to \s(2,k)\to 0$ mentioned in Theorem \ref{4.3}(3). Then $\Der L_1$ is 7-dimensional and has a basis  $(d_1,d_2,e,h,f,v_0,v_1)$ with multiplication
\begin{align*}
&[h,e]=-e+v_1,&&[h,v_0]=v_0,&& [e,v_0]=0,&&[f,v_0]=v_1,\\
&[h,f]=f,&&[h,v_1]=-v_1,&&[e,v_1]=v_0,&& [f,v_1]=0,\\
&[e,f]=h,&&[v_0,v_1]=0,
\end{align*}
and 
\begin{align*}
&d_1(e)=v_1,&&d_1(h)=d_1(f)=d_1(v_0)=d_1(v_1)=0,\\
&d_2(f)=v_0,&&d_2(h)=d_2(e)=d_2(v_0)=d_2(v_1)=0.
\end{align*}
\end{prop}
\begin{Proof}
Let $(e,h,f,v_0,v_1)$ denote the basis of $L_1$ mentioned in Theorem \ref{4.3}(3), and let $D$ denote any derivation of $L_1$. Note that $\rad(L_1)+D(\rad(L_1))$ is an ideal of $L_1$ of dimension $\le 2\dim\rad(L_1)=4<\dim L_1$. But $\rad(L_1)$ is the only proper ideal of $L_1$. Therefore $D(\rad(L_1))\subset \rad(L_1)$ holds. Then $D$ induces a derivation of $L_1/\rad(L_1)\cong \s(2,k)$. All derivations of this algebra are inner. We therefore may assume that 
$D(L_1)\subset \rad(L_1)$.
Put $D(h)=\alpha_0v_0+\alpha_1v_1$. Setting $D':=D+\ad(\alpha_0v_0-\alpha_1v_1)$ gives $D'(h)=0$. Considering $(\ad h)$-eigenvalues one obtains 
$$D'(e)=\beta_1v_1,\quad D'(f)=\beta_0v_0,\quad D'(v_i)=\gamma_iv_i\ (i=0,1).$$
Then 
\begin{align*}
0&=D'([h,e])-D'(-e+v_1)=[h,D'(e)]+D'(e)-\gamma_1v_1=-\gamma_1v_1,
\end{align*}
whence $\gamma_1=0$. Next,
$$\gamma_0v_0=D'(v_0)=D'([e,v_1])=[D'(e),v_1]+[e,D'(v_1)]=0,$$
hence $\gamma_0=0$. Then $D'=\beta_1d_1+\beta_0d_2$.
\\
 On the other hand, it is not hard to see that $d_1$ and $d_2$ as given in the proposition are derivations.
\end{Proof}

\begin{prop} \label{4.6}Let
$p=3$, and $L_1$ be the nonsplit extension $0\to V(1)\to L\to \s(2,k)\to 0$ mentioned in Theorem \ref{4.3}(3). There are up to algebra isomorphisms exactly 3  central extensions $$0\to k\to L\to  L_1\to 0,$$
namely every such algebra $L$ has a basis $(e,h,f,v_0,v_1,z)$ with multiplication
\begin{align*}
&[h,e]=-e+v_1,&&[h,v_0]=v_0,&&[e,v_0]=\alpha z, &&[f,v_0]=v_1,\\
&[h,f]=f,&&[h,v_1]=-v_1,&&[e,v_1]=v_0,&& [f,v_1]=0,\\
&[e,f]=h,&&[v_0,v_1]=\beta z,
\end{align*}
where one of the following occurs
$$ \alpha=\beta=0;\qquad  \alpha=0,\ \beta=1;\qquad
 \alpha=1,\ \beta=0.
$$
\end{prop}
\begin{Proof}
Let $kz$ denote the 1-dimensional center of $L$. By assumption, $L$ has a basis $(e,h,f,v_0,v_1,z)$ with multiplication
\begin{align*}
&[h,e]=-e+v_1,&&[h,v_0]=v_0,&& [e,v_0]=\alpha z,&&[f,v_0]=v_1,\\
&[h,f]=f,&&[h,v_1]=-v_1,&&[e,v_1]=v_0,&& [f,v_1]=\gamma z,\\
&[e,f]=h+\delta z,&&[v_0,v_1]=\beta z.
\end{align*}
Suppose $\beta\ne0$. Adjusting $z$ we may assume $\beta=1$. Substitute $e':=e+\alpha v_1$  to obtain $\alpha=0$.
\\
Suppose $\beta=0$. Then $\alpha=0$, or $\alpha\ne0$, in which case we set $z':=\alpha z$ to obtain $\alpha=1$.
\\
Finally, set $h':=h+\delta z$ to obtain $\delta=0$. We compute
$$0=[h,[e,f]]=[[h,e],f]+[e,[h,f]]=[-e+v_1,f]+[e,f]=[v_1,f]=-\gamma z.$$
This gives $\gamma=0$.

We show that the exposed algebras are nonisomorphic.  The radical of $L$ is nonabelian, if and only if $\beta=1$. Next suppose $\beta=0$. Then $L/\rad(L)$ acts on $\rad(L)$, and $\rad(L)$ is an indecomposable module if and only if   $\alpha\ne0$.
\end{Proof}

\section{Nonsolvable Lie algebras of dimension 6}

Dimension 6 is the lowest dimension for which parameter depending families of nonsolvable Lie algebras occur. Even more we face the fact that  the occurrence of nonrestricted 3-dimensional modules (for $p=3$)  gives rise to a rather long list of isomorphism types.

\begin{theo}\label{6.1} Let $L$ be a nonsolvable Lie algebra of dimension $6$ over a finite field $k$ of characteristic $2$.   One of the following occurs. 
\begin{enumerate}
\item $\dim \rad(L)=0:$
\begin{enumerate}
\item $L=W(1;\underline 2)^{(1)}\oplus W(1;\underline 2)^{(1)}$; 
\item $L=W(1;\underline 2)^{(1)}\otimes_k{\cal C}(L)$, where ${\cal C}(L)/k$ is a field extension of degree 2;
\end{enumerate}
\item $\dim \rad(L)=1:$ \\
$L=\Der W(1;\underline 2)^{(1)}\rtimes \rad(L)$, where $[W(1;\underline 2),\rad(L)]=\{0\}$ and $$\rad(L)=ku,\  [\partial^2,u]=\delta u,\ \delta\in\{0,1\};$$
\item $\dim \rad(L)=2:$ \\
$L=W(1;\underline 2)\rtimes \rad(L),$ where  $[W(1;\underline 2)^{(1)}, \rad(L)]=\{0\}$ and $$\rad(L)=kh\oplus ku,\quad  [h,u]=\delta u,\quad \delta\in\{0,1\},$$ and \begin{enumerate}
\item if $\delta=0:$ the action of $x^{(3)}\partial$ on $\rad(L)$ is given by one of the following matrices
\begin{align*} &\left(\begin{array}{cc}0&0\\ 0&0
\end{array}\right), \left(\begin{array}{cc}0&1\\ 0&0
\end{array}\right),\ \left(\begin{array}{cc}1&0\\ 0&1
\end{array}\right), \ \left(\begin{array}{cc}1&1\\ 0&1
\end{array}\right),\ 
  \left(\begin{array}{cc}0&\xi\\1&1
\end{array}\right)\ (\xi\in k^\ast).\end{align*}
\item if $\delta=1:$
 $$[x^{(3)}\partial,h]=0,\quad [x^{(3)}\partial,u]=\delta'u,\quad \delta'\in\{0,1\}.$$\end{enumerate}
\item $\dim \rad(L)=3:$  
\begin{enumerate}
\item $L=W(1;\underline 2)^{(1)}\oplus \rad(L)$ is the direct sum of two ideals, and 
$\rad(L)$ is given by Proposition \ref{3.5};
\item $L=W(1;\underline 2)^{(1)}\rtimes{\cal O}(1;\underline 2)/k$ is the semidirect sum of a subalgebra $\cong W(1;\underline 2)^{(1)}$ and the abelian ideal $\rad(L)\cong {\cal O}(1;\underline 2)/k$.
\item $L$ is the nonsplit extension
$0\to {\cal O}(1;\underline 2)/k\to L\to W(1;\underline 2)^{(1)}\to 0$; more exactly, $L$ has a basis $(e,h,f,v_1,v_2,v_3)$ with multiplication
\begin{align*}
&[h,e]=e+v_3,&&[h,v_1]=v_1,&&[e,v_1]=0,&&[f,v_1]=v_2,
\\& [h,f]=f,&&[h,v_2]=0,&&[e,v_2]=v_1,&&[f,v_2]=v_3,
\\& [e,f]=h,&&[h,v_3]=v_3,&&[e,v_3]=v_2,&&[f,v_3]=0,
\\&[v_i,v_j]=0. 
\end{align*}
\end{enumerate}
\end{enumerate}
The exposed algebras are mutually nonisomorphic.
 \end{theo}
 \begin{Proof}
(1) Consider the case $\rad(L)=\{0\}$. This case is covered by Theorems \ref{2.4} and \ref{3.1}. The following can occur.
\\[2mm] 
(a) $L=L_1\oplus L_2$, where $L_1,L_2$ are 3-dimensional simple. These algebras are isomorphic to $W(1;\underline 2)^{(1)}$.
\\[2mm]
(b) ${\cal C}(L)/k$ is a field extension of degree 2 and $L$ is 3-dimensional simple over ${\cal C}(L)$. Then $L$ is ${\cal C}(L)$-isomorphic to $W(1;\underline 2)^{(1)}$. Choose the ${\cal  C}(L)$-basis  $(e:=\partial,h:=x\partial,f:=x^{(2)}\partial)$ of $L$, which by definition  satisfies the equations
$$[e,h]=e,\ [e,f]=h,\ [h,f]=f.$$
Put $L':=ke\oplus kh\oplus kf$, which is a Lie algebra over $k$ isomorphic to $W(1;\underline 2)^{(1)}$. Then $L\cong L'\otimes_k{\cal  C}(L)$.
\\[2mm]
(2) Consider the case $\dim\rad(L)=1$. Let $\pi:L\to L/\rad(L)$ denote the canonical homomorphism. Then $\pi(L)$ is described in Theorem \ref{4.2}. As $\pi(L)$ is semisimple, only $\pi(L)=\Der W(1;\underline 2)^{(1)}$ is possible. This algebra contains a 3-dimensional ideal $Q$ isomorphic to $W(1;\underline 2)^{(1)}$. Set $L_0:=\pi^{-1}(Q)$, which is a 4-dimensional ideal of $L$ with 1-dimensional radical. Theorem \ref{4.1} yields 
$$L_0=P\oplus \rad(L),$$
where $P$ is an ideal of $L_0$ isomorphic to $W(1;\underline 2)^{(1)}$. There is a vector space decomposition
$$L=kd_1\oplus kd_2\oplus P\oplus ku,$$
where $P= L_0^{(1)}$ is an ideal of $L$, $ku=\rad(L)$, $[P,u]=\{0\}$, and (as $\pi(L)=\Der W(1;\underline 2)^{(1)}$) there is an isomorphism
$$\sigma: kd_1\oplus kd_2\oplus P\cong \Der W(1;\underline 2)^{(1)}.$$
We have to determine the action of $d_1:=\sigma^{-1}(\partial^2)$ and $d_2:=\sigma^{-1}(x^{(3)}\partial)$ on the 1-dimensional ideal $ku$. There are $r,s\in k$ not both 0 such that
$$[rd_1+sd_2,u]=0.$$
Note that $\sigma^{-1}(x\partial)\in P^{(1)}$ annihilates $u$.
Since $k$ is a perfect field of characteristic 2 one can choose $\alpha,\beta,\gamma,\delta\in k$ for which 
$$\gamma^2=r,\ \delta^2=s,\ \alpha\delta+\beta\gamma=1.$$
Proposition \ref{3.2} shows that there is an automorphism of $\Der W(1;\underline 2)^{(1)}$ which maps $x^{(3)}\partial$ onto $r\partial^2 +\gamma\delta x\partial+sx^{(3)}\partial$. Since $[\sigma^{-1}(r\partial^2+\gamma\delta x\partial+sx^{(3)}\partial),u]=0$, we may assume $[d_2,u]=0$.
\\
Next let $[d_1,u]=tu$ for some $t\in k$. If $t\ne0$, choose $\alpha\in k$ with $\alpha^2=t^{-1}$, set $\beta=\gamma=0$, $\delta=\alpha^{-1}$. With this choice there is an automorphism of $\Der W(1;\underline 2)^{(1)}$ which maps $\partial^2$ onto $t^{-1}\partial^2$ and $W(1;\underline 2)$ onto $W(1;\underline 2)$. Therefore we may assume $t=1$ in this case.
\\[2mm]
(3) Consider the case $\dim\rad(L)=2$. As before, let $\pi:L\to L/\rad(L)$ denote the canonical homomorphism.  As $\pi(L)$ is 4-dimensional and semisimple, only $\pi(L)= W(1;\underline 2)$ is possible (Theorem \ref{4.1}). This algebra contains a 3-dimensional ideal $Q$ isomorphic to $W(1;\underline 2)^{(1)}$. Set $L_0:=\pi^{-1}(Q)$, which is a 5-dimensional ideal of $L$ with 2-dimensional radical. Theorem \ref{4.2} yields that $L_0=P\oplus \rad(L)$ is the direct sum of ideals $P$ and $\rad(L)$, where $P\cong W(1;\underline 2)^{(1)}$ and 
\begin{align*}\rad(L)=kh\oplus ku, \ [h,u]=\delta u,\ \delta\in\{0,1\}.\end{align*}
 One has a vector space decomposition
$$L=kd\oplus P\oplus \rad(L),$$
and there is an isomorphism  $\sigma:kd\oplus P\cong W(1;\underline 2)$ (since  $\pi(L)= W(1;\underline 2)$). Note that $P=(kd\oplus P)^{(1)}\cong W(1;\underline 2)^{(1)}$ under this isomorphism.
We have to determine the action of $d:=\sigma^{-1}(x^{(3)}\partial)$ on $\rad(L)=kh\oplus ku$.
\\[2mm]
(a) Consider the case $\delta=0$: Let $\chi=T^2+\alpha T+\beta$ be the characteristic polynomial of $d$ acting on the abelian Lie algebra $ \rad(L)$. If  $\alpha\ne0$, we may adjust $x^{(3)}\partial$ by a nonzero scalar using Proposition \ref{3.2} as in former cases to obtain  $\alpha=1$. Similarly, if $d$ has a nonzero eigenvalue on $ \rad(L)$, then adjusting $x^{(3)}\partial$ by a nonzero scalar  we may assume that 1 is an eigenvalue.   This means that we may assume
$$\chi\in\{ T^2,\  T^2+1;\ T^2+T+\beta,\ (\beta\in k )\}.$$
If $\chi=T^2$, we may choose a basis $(u,v)$ of $ \rad(L)$ such that $[d,u]=[d,v]=0$ or choose a basis $(u,v)$ of $ \rad(L)$ such that  $[d,v]=u$, $[d,u]=0$ (describing the cases that $d$ acts semisimply or not on $\rad(L)$).
\\
If $\chi=T^2+1$, we similarly may  choose $u,v$ such that  $[d,u]=u$, $[d,v]=v$ or choose  $u,v$  such that  $[d,v]=u+v$, $[d,u]=u$.
\\
Suppose $\chi=T^2+T+\beta$. If all nonzero vectors are eigenvectors, then they are eigenvectors for the same eigenvalue $r$, and $\chi=(T-r)^2=T^2+r^2$ holds. This is not the requested form for $\chi$ in this case.
Therefore  there is a vector $u$ which is not an eigenvector. Then  we may  choose a basis $(u,v:=[d,u])$, and $[d,v]=(\ad d)^2(u)=\beta u+ v$ holds. Set $\xi:=\beta$.
\\[2mm]
(b) Consider the case $\delta=1$: Then $(\rad(L))^{(1)}=ku$ is an ideal of $L$. Set $[d,h]=\alpha h+\beta u$ and $[d,u]=\gamma u$ with $\alpha,\beta,\gamma\in k$. Then
\begin{align*}
\gamma u&=[d,u]=[d,[h,u]]=[\alpha h+\beta u,u]+[h,\gamma u]=(\alpha+\gamma)u.\end{align*}
Therefore $\alpha=0$. Substituting $d$ by $d+\beta u$ (these elements act identically on 
$W(1;\underline 2)^{(1)}$) one obtains $\beta=0$.   We may as in former cases adjust $x^{(3)}\partial$ to obtain $\gamma\in\{0,1\}$. 
 \\[2mm]
 (4) Consider the case $\dim\rad(L)=3$. Observe that $L/\rad(L)$ is 3-dimensional semisimple. Theorem \ref{3.1} yields $L/\rad(L)\cong W(1;\underline 2)^{(1)}$.
 \\[2mm]
 (a) Suppose $\rad(L)$ contains an ideal $I$ of $L$ different from $\{0\}$ and $\rad(L)$. Let $\pi:L\to L/I$ denote the canonical homomorphism.  Note that $\pi(L)/\rad \pi(L)\cong W(1;\underline 2)^{(1)}$.
 As $\dim\pi(L)\in\{4,5\}$, Theorems \ref{4.1} and \ref{4.2} show that $\pi(L)=Q\oplus \rad(\pi(L))$ is the direct sum of ideals where $Q\cong W(1;\underline 2)^{(1)}$.   Set $L_0:=\pi^{-1}(Q)$. Then $\dim L_0\in\{4,5\}$ and $L_0$ is an ideal of $L$ with radical $I$ of codimension 3. Theorems \ref{4.1} and \ref{4.2} yield 
$L_0=P\oplus I$, where $P\cong W(1;\underline 2)^{(1)}$ is an ideal of $L_0$. Since both $\rad(L)/I $, $I$ have dimension less than 3, Proposition \ref{3.2} yields that $P$ annihilates $\rad(L)$. Then $P$ is an ideal of $L$.
\\[2mm]
(b) Now suppose that $\rad(L)$ contains no ideal  of $L$ properly.  In particular, it is abelian, and $L/\rad(L)\cong W(1;\underline 2)^{(1)}$ acts on $\rad(L)$ irreducibly. The only faithful 3-dimensional module of $W(1;\underline 2)^{(1)}$ has been described in Proposition \ref{3.2}.  It is isomorphic to ${\cal O}(1;\underline 2)/k$. If the extension
$$0\to {\cal O}(1;\underline 2)/k\to L\to W(1;\underline 2)^{(1)}\to 0$$
 splits, then $L=W(1;\underline 2)^{(1)}\rtimes {\cal O}(1;\underline 2)/k$.
   \\[2mm]
 (c) Finally assume that $\rad(L)\cong {\cal O}(1;\underline 2)/k$ is abelian and the extension
$$0\to {\cal O}(1;\underline 2)/k\to L\to W(1;\underline 2)^{(1)}\to 0$$
 does not split. Let $\pi:L\to W(1;\underline 2)^{(1)}$ denote the canonical homomorphism. Choose  a basis  $(e,h,f,v_1,v_2,v_3)$ of $L$ such that
 $$\bar e:=\partial,\ \bar h:=x\partial,\ \bar f:=x^{(2)}\partial,\ v_1:=x+k,\ v_2:=x^{(2)}+k,\ v_3:= x^{(3)}+k,$$
 and $e,f$ are $(\ad h)$-root vectors.
These elements multiply as follows
\begin{align*}&[h,e]=e+\alpha_1v_1+\alpha_3v_3,&&[e,v_1]=0,&&[f,v_1]=v_2,&&[h,v_1]=v_1,\\
&[h,f]=f+\beta_1v_1+\beta_3v_3,&&[e,v_2]=v_1,&&[f,v_2]=v_3,&&[h,v_2]=0,\\
&[e,f]=h+\gamma_2v_2,&&[e,v_3]=v_2,&&[f,v_3]=0,&&[h,v_3]=v_3.
\end{align*}
Moreover,
\begin{align*}
0&=[h,[e,f]]=[e+\alpha_1v_1+\alpha_3v_3,f]+[e,f+\beta_1v_1+\beta_3v_3]\\&=h+\alpha_1v_2+h+\beta_3v_2=(\alpha_1+\beta_3)v_2,
\end{align*}
whence
\begin{align}
\alpha_1=\beta_3.
\end{align}
Next let $r,s\in k$ and set
\begin{align*}
&e':=e,&&v_1':=v_1,
\\&h':=h+rv_2,&&v_2':=v_2,\\&f':=f+se,&&v_3':=v_3+sv_1.
\end{align*}
It is not hard to compute
\begin{align*}&[e',v'_1]=0,&&[f',v'_1]=v'_2,&&[h',v'_1]=v'_1,\\
&[e',v'_2]=v'_1,&&[f',v'_2]=v'_3,&&[h',v'_2]=0,\\
&[e',v'_3]=v'_2,&&[f',v'_3]=0,&&[h',v'_3]=v'_3,
\end{align*}
and $(\ad h'-{\rm Id})^2(e')=(\ad h'-{\rm Id})^2(f')=0$. In addition, 
\begin{align*}
[h',f']&=[h,f+se]+[rv_2,f+se]\\
&=(f+\beta_1v_1+\beta_3v_3)+s(e+\alpha_1v_1+\alpha_3v_3)+rv_3+srv_1
\\&=f+se+(s\alpha_1+\beta_1+sr)v_1+(s\alpha_3+\beta_3+r)v_3
\\&=f'+((s\alpha_1+\beta_1+sr)+s(s\alpha_3+\beta_3+r))v'_1+(s\alpha_3+\beta_3+r)v'_3.
\end{align*}
If $\alpha_3\ne0$, then solve the equations $$s^2=\alpha_3^{-1}\beta_1,\ r=s\alpha_3+\beta_3.$$
One gets $s\alpha_3+\beta_3+r=0$ and because of equation (5.1)
$$s\alpha_1+\beta_1+sr=s\alpha_1+\beta_1+s(s\alpha_3+\beta_3)=\beta_1+s^2\alpha_3=0.$$
Therefore we obtain in this case $[h',f']=f'$. \\
If $\alpha_3=0$, set $r=\alpha_1$ and $s=0$. Then $[h',e']=e'$ holds, and interchanging $f'$ with $e'$ and $v_1'$ with $v_3'$ again gives $[h',f']=f'$. Using equation (5.1) one obtains in both cases a multiplication
$$[h',e']=e'+\alpha v_3',\quad [h',f']=f',\quad [e',f']=h'+\gamma v_2'.$$
Next, set $f'':=f'+\gamma v_3'$, $h'':=h'$, $e'':=e'$,  and obtain $[e'',f'']=h''$. Since the sequence does not split by assumption, one has $\alpha\ne0$. Now set
$$v_1'':=\alpha v_1',\ v_2'':=\alpha v_2',\ v_3'':=\alpha v_3'$$
to obtain $\alpha=1$.
\\[2mm]
(5) It remains to show that the listed algebras are nonisomorphic. 
 It is clear that we only have to discuss algebras within the same number 1 - 4.
\\[2mm]
(i) The algebra of 1.(a) is nonsimple, while that of 1.(b) is so.
\\[2mm]
(ii) The algebras in 2. have a unique minimal simple ideal $W(1;\underline 2)^{(1)}$, and the quotient $L/W(1;\underline 2)^{(1)}$ is abelian if and only if $\delta=0$.
\\[2mm]
(iii) The algebras in 3. have a unique minimal simple ideal $ W(1;\underline 2)^{(1)}$, and $L/W(1;\underline 2)^{(1)}$ is 3-dimensional solvable. In addition, they have abelian radical if and only if $\delta=0$. Apply the nonisomorphism claim of Proposition \ref{3.5} for the cases $\delta=0$ and $\delta=1$ separately.
\\[2mm]
(iv) The algebras mentioned in 4.(a) and 4.(b) are obviously not isomorphic. It remains to show that the algebra mentioned in (c) in fact is nonsplit. In order to do so we assume that $L$ has a basis and multiplication as listed but contains a subalgebra $P\cong W(1;\underline 2)^{(1)}$. Then $L=P\rtimes \rad(L)$. Let $(e',h',f')$ denote the projections of $(e,h,f)$ into $P$. The multiplication of $L/\rad(L)$ shows that 
$$[h',e']=e',\ [h',f']=f',\ [e',f']=h'.$$
Write
$$e=e'+\alpha_1v_1+\alpha_2v_2+\alpha_3v_3,\ h=h'+\beta_1v_1+\beta_2v_2+\beta_3v_3,\ f=f'+\gamma_1v_1+\gamma_2v_2+\gamma_3v_3.$$
Then
\begin{align*}
e+v_3&=[h,e]=[h'+\beta_1v_1+\beta_2v_2+\beta_3v_3,e'+\alpha_1v_1+\alpha_2v_2+\alpha_3v_3]\\&=[h',e']+[\beta_1v_1+\beta_2v_2+\beta_3v_3,e']+[h',\alpha_1v_1+\alpha_2v_2+\alpha_3v_3]\\
&=e'+\beta_2v_1+\beta_3v_2+\alpha_1v_1+\alpha_3v_3\\&=e+\alpha_2v_2+\beta_2v_1+\beta_3v_2,
\end{align*}
 a contradiction.
\end{Proof}

When $char(k)>2$, we have to split the investigation into several cases.

 \begin{theo}\label{6.2}
Let $L$ be a nonsolvable Lie algebra of dimension $6$ over a finite field $k$ of characteristic $p\ge3$.   Assume $\dim\rad(L)\le2$. Then  $\dim\rad(L)\le1$ and one of the following occurs.
\begin{enumerate}
\item $\dim \rad L=0:$ 
\begin{enumerate}
\item $L=\s(2,k)\oplus\s(2,k)$;
\item $L=\s(2,{\cal C}(L))$, where ${\cal C}(L)/k$ is a field extension of degree 2;
\end{enumerate}
\item $\dim \rad L=1:$ $p=5$ and 
\begin{enumerate}
\item $L=W(1;\underline 1)\oplus C(L)$;
\item $L$ is the nonsplit central extension of $W(1;\underline 1)$, i.e., $$L=\Big(\oplus_{i=-1}^3ke_i\Big)\oplus kz \text{ where }[L,z]=\{0\}$$ and
$$ [e_i,e_j]=\left\{\begin{array}{ll}(j-i)e_{i+j}&\text{ if }-1\le i+j\le3,\\ z&\text{ if }i=2,j=3,\\-z&\text{ if }i=3,j=2,
\\ 0&\text{ otherwise}.\end{array}\right.$$
\end{enumerate} 
\end{enumerate}
The exposed algebras are mutually nonisomorphic.
\end{theo}
\begin{Proof}
(1) The case $\rad(L)=\{0\}$  is  covered by Theorems \ref{2.4} and \ref{3.1}. The following can occur.
\\[2mm] 
(a) $L=L_1\oplus L_2$, where $L_1,L_2$ are 3-dimensional simple. These algebras are isomorphic to $\s(2,k)$.
\\[2mm]
(b) ${\cal C}(L)/k$ is a field extension of degree 2 and $L$ is 3-dimensional simple over ${\cal C}(L)$. Then $L$ is ${\cal C}(L)$-isomorphic to $\s(2,{\cal C}(L))$.\\[2mm]
(2) Suppose $\dim \rad(L)=1$. Then $L/\rad(L)$ is 5-dimensional and semisimple. Theorem \ref{4.2} shows that $p=5$ and $L/\rad(L)\cong W(1;\underline 1)$. The central extensions of the Witt algebra can be found for example in \cite[page 428]{St04} (the proof there does not at all need an algebraically closed ground field). They  are as claimed.
\\[2mm]
(3) Suppose $\dim \rad(L)=2$. Then $L/\rad(L)$ is 4-dimensional and semisimple. However, Theorem \ref{4.1} shows that no such algebra exists.
 \\[2mm]
(4) Obviously, there are no isomorphisms among the exposed algebras.

\end{Proof}

If $\rad(L)$ is at least 3-dimensional, a great variety of algebras does occur (mainly for $p=3$). Since no Lie algebra of dimension less than 3 is semisimple it can only be that $\dim \rad(L)=3$ in the present case.

\begin{theo}\label{6.3}
Let $L$ be a nonsolvable Lie algebra of dimension $6$ over a finite field $k$ of characteristic $p\ge3$.   Assume that $\dim\rad(L)= 3$ and the extension
$$0\to \rad(L)\to L\to L/\rad(L)\to 0$$
splits.  Then $L=P\rtimes\rad(L)$ is the semidirect sum of a subalgebra $P\cong \s(2,k)$ and the ideal $\rad(L)$, and one of the following occurs.
\begin{enumerate}
\item $L=P\oplus \rad(L)$ is the direct sum of two ideals, and 
$\rad(L)$ is given by Proposition \ref{3.5};
\item  $\rad(L)$ is abelian and completely reducible as a $P$-module of the form 
$$V(0)\oplus V(1)\quad\text{or}\quad V(2);$$
\item $\rad(L)$ is abelian, $p=3$ and  
\begin{enumerate}
\item  $\rad(L)=V(2,\chi)$ for some $\chi\in P^\ast$ described in Proposition \ref{3.4}, or
\item $L=W(1;\underline 1)\rtimes {\cal O}(1;\underline 1)$, or
\item $L=W(1;\underline 1)\rtimes {\cal O}(1;\underline 1)^\ast$;
\end{enumerate}
\item  $\rad(L)=V(1)\oplus C(L)$ decomposes as a $P$-module and $(\rad(L))^{(1)}=C(L)$ is 1-dimensional;
\item $p=3$, $\rad(L)\cong {\cal O}(1;\underline 1)$ as a $P$-module, and $[x,x^{(2)}]=1$, $[L,1]=\{0\}$;
\item
$\rad(L)=kd\oplus (\rad(L))^{(1)}$, $(\rad(L))^{(1)}$ is 2-dimensional abelian and isomorphic to $V(1)$ as a $P$-module, and $[P,d]=\{0\}$, $ [d,v]=v\ \forall v\in \rad(L)$;
\end{enumerate}
All algebras listed are mutually nonisomorphic.
\end{theo}
\begin{Proof}  The present assumption means that $L$ has a subalgebra $P\cong L/\rad(L)$. This subalgebra is semisimple and at most 3-dimensional. Therefore Theorem \ref{3.1} shows that $P\cong \s(2,k)$. Moreover, $\rad(L)$ is 3-dimensional, and therefore is determined by Proposition \ref{3.5}. We have to describe the action of $P$ on $\rad(L)$.\\
The case $[P,\rad(L)]=\{0\}$ is listed as 1. of this theorem. \\
Otherwise the simplicity of $P$ ensures that $P$ acts faithfully on $\rad(L)$.\\[2mm]
(1) Suppose $\rad(L)$ is abelian. \\
(i) If the $P$-module $\rad(L)$ is completely reducible, then the following is possible
$$\rad(L)=V(0)\oplus V(0)\oplus V(0), \ V(0)\oplus V(1),\ \text{irreducible}.$$
In the first case $P$ annihilates $\rad(L)$. But we assume that $P$ acts faithfully on $\rad(L)$. 
The second case is listed as 2. of the theorem.
\\[2mm]
Consider the case that $\rad(L)$ is $P$-irreducible.    Recall that  ${\cal C}(\rad(L),P)$ is an extension field of $k$, and $\rad(L)$ is a vector space over ${\cal C}(\rad(L),P)$. Let $d_1$ denote the $k$-dimension of ${\cal C}(\rad(L),P)$ and $d_2$ denote the ${\cal C}(\rad(L),P)$-dimension of $\rad(L)$. Then $d_1d_2=\dim_k\rad(L)=3$. The case   $d_1>1$ gives $d_2=1$, whence $P$ would act abelian on $\rad(L)$, a contradiction. Hence  $d_1=1$, which means ${\cal C}(\rad(L),P)=k$. In the course of section 3 we have mentioned that the $P$-module $\rad(L)$ is one of the following. Namely, if $p>3$, then it is of type $V(2)$, while for $p=3$ it is of type $V(2,\chi)$ for some $\chi\in P^\ast$. The ${\rm Aut}\,\s(2,k)$-orbits of nonzero irreducible characters are ruled by Proposition \ref{3.4}. These cases are listed in 2. and 3.(a).
\\[2mm]
(ii) If the $P$-module $\rad(L)$ is decomposable, then $\rad(L)=U_1\oplus U_2$, where $U_1,U_2$ are $P$-modules and one of these is 1-dimensional while the other is 2-dimensional. But 2-dimensional modules are completely reducible,  so this is the former case (i).
\\[2mm]
(iii) Suppose the $P$-module $\rad(L)$ is indecomposable and not irreducible. Then there is a submodule $U$ of dimension 1 or 2 and the quotient module has dimension 2 or 1, respectively. Therefore passing to the dual module we may assume that $\dim U=1$ and $\rad(L)/U\cong V(1)$. 
\\
Set $U=kz$. Take an $\s(2)$-triple $(e,h,f)$ of $P$, choose accordingly a basis $(\bar v_0,\bar v_1)$ of $\rad(L)/U$ and take inverse images $v_0,v_1$ of these as eigenvectors for $\ad h$. Then 
 $(v_0,v_1,z)$ is a basis of $\rad(L)$ and the following equations hold
 \begin{align*}&[h,e]=2e,&&[h,v_0]=v_0,&&[e,v_0]=\alpha z,&&[f,v_0]=v_1,\\
 &[h,f]=-2f,&&[h,v_1]=-v_1,&& [e,v_1]=v_0,&&[f,v_1]=\beta z,\\
 &[e,f]=h,&&  [v_0,v_1]=0,&&[L,z]=\{0\}.\end{align*}
We have assumed that the module does not split, so we assume that $\alpha\ne0$ or $\beta\ne0$. Considering $(\ad h)$-eigenvalues one obtains $p=3$. Intertwining $e,f$ and $v_0,v_1$ if necessary gives $\alpha\ne0$ in any case. Adjusting $v_0,v_1$ by $\alpha$ gives $\alpha=1$. Since $k$ is perfect, there is $\kappa\in k$ satisfying $\kappa^3=\beta$. Put
\begin{align*}
&e':=e,&&h':=h+\kappa e,&&f':=f+\kappa h-\kappa^2e,\\&
v_0':=v_0+\kappa z,&&v_1':=v_1+\kappa v_0-\kappa ^2z,&&z':=z.
\end{align*}
Then (as $p=3$)
\begin{align*}
&[h',e']=2e', \\
&[h',f']=-2f-2\kappa ^2e+\kappa h-2\kappa ^2e=-2(f+\kappa h-\kappa ^2e)=-2f',\\
&[e',f']=h',\\
&[e',v_0']=z', \\
&[h',v_0']=v_0+\kappa z=v_0',\\
&[f',v_0']=v_1+\kappa v_0-\kappa ^2z=v_1',\\
&[e',v_1']=v_0+\kappa z=v_0',\\
& [h',v_1']=(-v_1+\kappa v_0)+\kappa(v_0+\kappa z)=-v_1',\\
&[f',v_1']=(\beta z+\kappa v_1)+\kappa (-v_1+\kappa v_0)-\kappa ^2(v_0+\kappa z)=(\beta-\kappa ^3)z=0.
\end{align*}
This brings us to the case $\beta=0$. It is now not hard to see that the linear mapping
\begin{align*}
&\tau:L\to W(1;\underline 1)\rtimes {\cal O}(1;\underline 1),\\
&\tau(z):=1,\ \tau(v_0):=x,\ \tau(v_1):=x^{(2)},\\
&\tau(e):=\partial,\ \tau(h):=x\partial,\ \tau(f):=x^{(2)}\partial
\end{align*}
is in fact a Lie algebra isomorphism. This case is listed in 3.(b), and the dual case is listed in 3.(c).
\\[2mm]
(2) Suppose $\rad(L)$ is as in case 2 of Proposition \ref{3.5}. Then $\rad(L)$ has a $P$-composition series $\{0\}\subset C(\rad(L))\subset C(\rad(L))+(\rad(L))^{(1)}\subset \rad(L)$ with 1-dimensional factors.  Then $P$ annihilates all these, and therefore $P$ annihilates $\rad(L)$. But this is not true under our assumptions.
\\[2mm]
(3) Suppose $\rad(L)$ is as in case 3 of Proposition \ref{3.5}. Then $\rad(L)$ is not $P$-irreducible, but contains a 1-dimensional submodule. With this modification the proof of (1) applies verbally, and gives the cases 4 and 5 of the theorem.
\\[2mm]
(4) Suppose $\rad(L)=kd\oplus (\rad(L))^{(1)}$ is as in case 4 of Proposition \ref{3.5}. Then  $d$ acts invertibly on $(\rad(L))^{(1)}$ and $(\rad(L))^{(1)}$ is 2-dimensional abelian. If $(\rad(L))^{(1)}$ is not $P$-irreducible, then $\rad(L)$ has a $P$-composition series with only 1-dimensional factors. In this case $P$ annihilates the factors and hence annihilates $\rad(L)$. But this does not happen under the general assumption. Therefore $(\rad(L))^{(1)}$ is $P$-irreducible. In addition, $\rad(L)/(\rad(L))^{(1)}$ is 1-dimensional, hence it is annihilated by $P$. Therefore we have that $[d,P]\subset (\rad(L))^{(1)}$.  Decompose $L$ into the Fitting components with respect to $\ad d$. As a result of the former deliberations, $L_1(\ad d)=(\rad(L))^{(1)}$, $L_0(\ad d)\cong L/(\rad(L))^{(1)}\cong kd\oplus P$. This now means that $L=Q\rtimes (\rad(L))^{(1)}$ is the semidirect sum of a subalgebra $Q\cong kd\oplus P$ and the ideal $(\rad(L))^{(1)}$. As a consequence, $[d,P]=\{0\}$.
\\
Recall that  ${\cal C}((\rad(L))^{(1)},P)$ is an extension field of $k$, and $(\rad(L))^{(1)}$ is a vector space over ${\cal C}((\rad(L))^{(1)},P)$. Let $d_1$ denote the ${\cal C}((\rad(L))^{(1)},P)$-dimension of $\rad(L)$ and $d_2$ denote the $k$-dimension of ${\cal C}((\rad(L))^{(1)},P)$. Then $$d_1d_2=\dim_k(\rad(L))^{(1)}=2.$$ The case   $d_2>1$ gives $d_1=1$, whence $P$ would act abelian on $(\rad(L))^{(1)}$, a contradiction. Hence  $d_2=1$, which means ${\cal C}((\rad(L))^{(1)},P)=k{\rm Id}$. Note that $\ad d$ is contained in ${\cal C}((\rad(L))^{(1)},P)$. Then one can adjust $d$ to obtain that $d$ acts as the identity on $(\rad(L))^{(1)}$. This is case 6 of the theorem.
\\[2mm]
(5) The set $\{x\in L\mid [x,\rad(L)]=\{0\}\}$ is nonsolvable in case 1 and solvable in all other cases of the theorem. This observation and 
properties of the radical distinguish all algebras listed in the theorem.

\end{Proof}

\begin{theo}\label{6.4}
Let $L$ be a nonsolvable Lie algebra of dimension $6$ over a finite field $k$ of characteristic $p\ge3$.   Assume that $\dim\rad(L)= 3$ and the extension
$$0\to \rad(L)\to L\to L/\rad(L)\to 0$$
does not split.  Then $p=3$, and $L$ is one of the algebras described in Proposition \ref{4.6} or $L$ has a 5-dimensional ideal $L_1$ and a basis  $(d ,e,h,f,v_0,v_1)$ with multiplication
\begin{align*}
&[h,e]=-e+v_1,&&[h,v_0]=v_0,&& [e,v_0]=0,&&[f,v_0]=v_1,\\
&[h,f]=f,&&[h,v_1]=-v_1,&&[e,v_1]=v_0,&& [f,v_1]=0,\\
&[e,f]=h,&&[v_0,v_1]=0,
\end{align*}
and 
\begin{align*}
&[d,e]=v_1,&&[d,h]=[d,f]=[d,v_0]=[d,v_1]=0,\\
\intertext{or}
&[d,f]=v_0,&&[d,h]=[d,e]=[d,v_0]=[d,v_1]=0.
\end{align*}
These algebras are mutually nonisomorphic.
\end{theo}
\begin{Proof}The nonsplitting of the extension means that $L$ has no subalgebra isomorphic to $\s(2,k)$.\\[2mm]
(1) Suppose $\rad(L)$ has an $L$-composition series with only 1-dimensional factors. Let $I$ denote a 1-dimensional ideal and $\pi:L\to L/I$ be the canonical homomorphism. Then $\pi(L)$ is 5-dimensional with 2-dimensional radical $\rad(L)/I$, and this radical is not $\pi(L)$-irreducible. Therefore $\pi(L)$ is described in case 2 of  Theorem \ref{4.3}. This theorem shows that $\pi(L)$ contains a subalgebra $P'\cong\s(2,k)$. Then $\pi^{-1}(P')$ is 4-dimensional and isomorphic to $\frak{gl}(2,k)$ (Theorem \ref{4.1}). As a consequence, $L$ contains a subalgebra $P\cong \s(2,k)$. But then the extension splits, a contradiction.  
\\[2mm]
(2) Suppose $\rad(L)\supset I\supset \{0\}$ is an $L$-composition series with  $\dim I=1$, $\dim  \rad(L)/I=2$. Note that $L/I$ has an irreducible 2-dimensional radical, hence is ruled by cases 2(c) or 3 of Theorem \ref{4.3}. In case 2(c) $L/I$ contains a subalgebra $Q\cong\s(2,k)$. The inverse image of $Q$ in $L$ is 4-dimensional. Theorem \ref{4.1} shows that $L$ contains a subalgebra isomorphic to $\s(2,k)$, a contradiction. Therefore case 3 of Theorem \ref{4.3} applies. Then $L$ is a central extension of $L/I$, and therefore is described in Proposition \ref{4.6}.
\\[2mm]
 (3)  Suppose  $\rad(L)\supset I\supset \{0\}$ is an $L$-composition series with  $\dim I=2$, $ \dim \rad(L)/I=1$. Note that $I$ is abelian as it is $L$-irreducible. Theorem \ref{4.1} shows that $L/I\cong \frak{gl}(2,k)$. 
The inverse image $L_1$ of $\s(2,k)$ in $L$ is a 5-dimensional ideal of $L$. 
In addition,   $\rad(L)$ contains an element $d$ so that the following holds $$L=kd\oplus L_1,\ \rad(L)=kd\oplus I,\ [d,L]\subset I.$$ 
(i) Suppose $[d,I]\ne\{0\}$. Then $[d,I]=I$ because $I$ is $L$-irreducible. Decompose $L=L_0(\ad d)\oplus L_1(\ad d)$ into its Fitting components with respect to $\ad d$.  Obviously, $L_1(\ad d)=I$ and therefore $L_0(\ad d)\cong L/I\cong \frak{gl}(2,k)$. But then there is a subalgebra $P\subset L_0(\ad d)$ isomorphic to $\s(2,k)$, a contradiction. 
 \\[2mm]
 (ii) Suppose $[d,I]=\{0\}$. If $[d,L]=\{0\}$, then $L$ has the 1-dimensional ideal $kd$. This case has been treated in (2).\\
Therefore we assume $[d,L]\ne\{0\}$. Since $\rad(L_1)=I$ is $L_1$-irreducible, $L_1$ is described by cases 2(c) or 3 of Theorem \ref{4.3}. But in case 2(c) the extension would split. Therefore we are in case 3, and $L$ is a 6-dimensional subalgebra of $\Der L_1$. The latter algebra has been described in Proposition \ref{4.4}. According to that proposition $L$ has a basis $(d ,e,h,f,v_0,v_1)$ with multiplication
\begin{align*}
&[h,e]=-e+v_1,&&[h,v_0]=v_0,&& [e,v_0]=0,&&[f,v_0]=v_1,\\
&[h,f]=f,&&[h,v_1]=-v_1,&&[e,v_1]=v_0,&& [f,v_1]=0,\\
&[e,f]=h,&&[v_0,v_1]=0,
\end{align*}
and $d=\alpha_1 d_1+\alpha_2d_2$, whence
\begin{align*}
&[d,e]=\alpha_1 v_1,&&[d,f]=\alpha_2 v_0,&&d(h)=d(v_0)=d(v_1)=0.
\end{align*}
If $\alpha_1=0$ or $\alpha_2=0$, then the algebra is listed in the theorem. Therefore we assume $\alpha_1\alpha_2\ne0$.
Solve the equation 
$$q^3=\alpha_1^{-1}\alpha_2,\quad \text{  set } s:=q^2,\ \beta:=-\alpha_1\alpha_2^{-1}s.$$
Then
$$s^3=\alpha_1^{-2}\alpha_2^2,\quad \beta^2s=1,\quad \alpha_2s^{-1}\beta=-\alpha_1$$
hold. Note that $L_1=ke\oplus \ker(\ad f)^2$ and $(\ad f)^2(ke)=kf$. For $x\in L_1$ and $y\in\ker(\ad f)^2$ one computes
$$[(\ad f)^2(x),(\ad f)(y)]=[(\ad f)(x),(\ad f)^2(y)]=0.$$
 Therefore $\exp(\ad \beta f)$ is an automorphism of $L_1$. Set
\begin{align*}
&e':=\exp(\ad \beta f)(e)=e-\beta h-\beta^2 f,&&v_0':=\exp(\ad \beta f)(v_0)=v_0+\beta v_1,\\
&h':=\exp(\ad \beta f)(h)=h-\beta f,&&v_1':=\exp(\ad \beta f)(v_1)=v_1,\\
&f':=\exp(\ad \beta f)(f)=f,
\end{align*}
and $$D:=\alpha_1d_1+\alpha_2d_2+\alpha_1sv_0+\alpha_2s^{-1}v_1\equiv d\pmod{L_1}.$$
Then 
\begin{align*}
[D,e']&=\alpha_1v_1-\alpha_2\beta^2v_0+\alpha_1s(\beta v_0+\beta^2v_1)+\alpha_2s^{-1}(-v_0-\beta v_1)\\
&=(-\alpha_2\beta^2+\alpha_1s\beta-\alpha_2s^{-1})v_0+(\alpha_1+\alpha_1s\beta^2-\alpha_2s^{-1}\beta)v_1\\
&=s^{-1}(-\alpha_2-\alpha_1^2\alpha_2^{-1}s^3-\alpha_2)v_0+3\alpha_1v_1=0,\\
[D,h']&=-\alpha_2\beta v_0+\alpha_1s(-v_0+\beta v_1)+\alpha_2s^{-1}v_1\\&=(\alpha_1s-\alpha_1s)v_0+s^{-1}(-\alpha_1^2\alpha_2^{-1}s^3+\alpha_2)v_1=0,\\
[D,f']&=\alpha_2v_0-\alpha_1sv_1=\alpha_2(v_0-\alpha_1\alpha_2^{-1}sv_1)=\alpha_2v_0',\\
[D,v_0']&=[D,v_1']=0.
\end{align*}
Therefore the linear mapping $kd\oplus L_1\to kd_2\oplus L_1$ given by $\exp(\ad \beta f)$ on $L_1$ and $D\mapsto \alpha_2d_2$ is an algebra isomorphism.
\\[2mm]
 (4) Suppose $\rad(L)$ is $L$-irreducible. Then it is abelian and is an irreducible $L/\rad(L)$-module. In the course of section 3 we have mentioned that it is one of the following. Namely, if $p>3$, then it is of type $V(2)$, while for $p=3$ it is of type $V(2,\chi)$ for some linear form $\chi:\s(2,k)\to \bar k$. 
 
 Consider first the case that $p=3$ and $\chi\ne0$. Observe that $\bar L:=L/\rad(L)\cong \s(2,k)$ is a restricted Lie algebra. There is $x\in L\setminus \rad(L)$ such that 
 $$(\ad \bar x\big|_{\rad(L)})^3-\ad \bar x^{[3]}\big|_{\rad(L)}=\chi(\bar x)^3{\rm Id}_{\rad(L)}\ne0.$$
 Choose $y\in L$ for which $\bar y=\bar x^{[3]}$, set $D:=(\ad x)^3-\ad y\in \Der L$. Then 
 $$D(L)\subset \rad(L),\ D(v)=\chi(\bar x)^3v\ne 0\quad\forall v\in\rad(L).$$
 Decompose $L=L_0(D)\oplus L_1(D)$ into the Fitting components with respect to $D$. The above means $L_1(D)=\rad(L)$. Then $L_0(D)\cong L/\rad(L)$ is a subalgebra isomorphic to $\s(2,k)$. But then the extension in question splits.
 \\
 As a consequence, $\rad(L)\cong V(2)$ as an $L/\rad(L)$-module in all cases.  Choose an $\s(2)$-triple $(\bar e,\bar h,\bar f)$ in $L/\rad(L)$ and preimages $e,h,f$ such that $e,f$ are $h$-root vectors. Next choose 
 a basis $(v_0,v_1,v_2)$ of $\rad(L)$ so that 
\begin{align*}
&[ h,v_0]=2v_0,&&[ e,v_0]=0,&&[ f,v_0]=v_1,\\&[ h,v_1]=0,&&[ e,v_1]=2v_0,&&[ f,v_1]=v_2,\\&[ h,v_2]=-2v_2,&&[ e,v_2]=2v_1,&&[f,v_2]=0.
\end{align*}
Since $\{e,v_0\}$, $\{h,v_1\}$ and $\{f,v_2\}$ span the respective  $h$-root spaces, one has
$$[h,e]=2e+\alpha v_0,\quad [h,f]=-2f+\beta v_2,\quad[e,f]=h+\gamma v_1$$
for some $\alpha,\beta,\gamma\in k$. Note that
\begin{align*}
0&=[h,[e,f]]=[2e+\alpha v_0,f]+[e,-2f+\beta v_2]=-\alpha v_1+2\beta v_1.
\end{align*}
This gives $\alpha=2\beta$. 
Set
$$e':=e+(\gamma-\beta)v_0,\quad  h':=h+\beta v_1,\quad f':=f.$$
Then 
\begin{align*}
&[h',e']=(2e+\alpha v_0)+2(\gamma-\beta)v_0-2\beta v_0=2e',\\
&[h',f']=(-2f+\beta v_2)-\beta v_2=-2f',\\
&[e',f']=(h+\gamma v_1)-(\gamma-\beta)v_1=h'.
\end{align*}
Hence $P:=ke'+kh'+kf'$ is a subalgebra of $L$ isomorphic to $\s(2,k)$. Therefore this case does not occur. 
\\[2mm]
(5) We have to show that the exposed algebras are nonisomorphic. The algebras of Proposition \ref{4.6} have a 1-dimensional center, while the other algebras  listed in the theorem are centerless. It has been proved that the 3 types of algebras of Proposition \ref{4.6} are nonisomorphic. We show that the  algebras $kd_1\oplus L_1$ and $ kd_2\oplus L_1$ are nonisomorphic. Suppose on the contrary that there exists an algebra isomorphism $\sigma:kd_1\oplus L_1\to kd_2\oplus L_1$. Since  $(kd_1\oplus L_1)^{(1)}=(kd_2\oplus L_1)^{(1)}=L_1$, 
$\sigma$ induces by restriction an automorphism of $L_1$, and therefore it maps the unique minimal ideal $\rad(L_1)$ of $L_1$ onto itself. Set $$\sigma(d_1):=D=\alpha d_2+\beta_0 e+\beta_1 h+\beta_2 f+\gamma_0v_0+\gamma_1v_1$$
  and $${\cal A}:=\{x\in L_1\mid [d_1,x]=0\}=kh+kf+kv_0+kv_1.$$
 Then $\sigma({\cal A})=\{x\in L_1\mid [D,x]=0\}$. Since $\rad(L_1)\subset {\cal A}$, the above reasoning shows that $\rad(L_1)\subset \sigma({\cal A})$, and this gives $\beta_0=\beta_1=\beta_2=0$. Moreover, $\sigma({\cal A})/\rad(L)$ is 2-dimensional, and therefore there exists an element $x\in \sigma({\cal A})$ of the form
$x=\delta_0e+\delta_1 h\ne0$.
We compute
\begin{align*}
0&=[D,x]=[\alpha d_2+\gamma_0 v_0 +\gamma_1v_1,\delta_0e+\delta_1h]=(-\gamma_0\delta_1-\gamma_1\delta_0-\gamma_1\delta_1)v_0+\gamma_1\delta_1v_1,
\end{align*} whence
$$
\gamma_0\delta_1+\gamma_1\delta_0=0,\quad
\gamma_1\delta_1=0.$$
Both cases $\delta_1=0$ and $\delta_1\ne0$ yield $\gamma_1=0$. Consequently, $D=\alpha d_2+\gamma_0v_0$ and $\sigma({\cal A})\subset ke+k(\alpha h-\gamma_0f)+kv_0+kv_1$. But since ${\cal A}$ is 4-dimensional, this is only possible if $\gamma_0=0$,  
$$D=\alpha d_2,\qquad \sigma({\cal A})=ke+kh+kv_0+kv_1.$$
Next we oberve that $f\in {\cal A}^{(1)}\setminus \rad(L_1)$, and therefore $\sigma(f)\in \sigma({\cal A})^{(1)}\setminus \rad(L_1)=(ke+kv_0+kv_1)\setminus kv_0+kv_1$. Arguing similarly for $h$ we obtain
\begin{align*}
&\sigma(f)=\lambda e+\mu_0v_0+\mu_1v_1,\quad \lambda\ne0,\\
&\sigma(h)=\kappa h+\lambda' e+\mu_0'v_0+\mu_1'v_1,\quad \kappa\ne0.
\end{align*}
Then 
\begin{align*}
\lambda e+\mu_0v_0+\mu_1v_1&=
\sigma(f)=\sigma([h,f])=[\sigma(h),\sigma(f)]\\&=[\kappa h+\lambda' e+\mu_0'v_0+\mu_1'v_1,\lambda e+\mu_0v_0+\mu_1v_1]\\
&\equiv \kappa\lambda(-e+v_1)-\kappa\mu_1v_1\pmod{kv_0},
\end{align*}
and this gives $\kappa=-1$, $\mu_1=-\lambda+\mu_1$. But then $\lambda=0$, a contradiction.
\end{Proof}

Counting the isomorphism classes of this section we obtain the number of isomorphism classes of nonsolvable 6-dimensional Lie algebras over a finite field $k$ as follows
\begin{align*}
&p=2:&&15+2|k|,
\\
&p=3:&&19+|k|+\big|\{\xi\in k\mid T^3+T^2=\xi\text{ has a solution in }k\}\big|,\\
&p=5:&&12+|k|,\\
&p>5:&&11+|k|.
\end{align*}


\begin{thebibliography}{99}


\bibitem{Bl} R.E. Block,
Determination of the differentiably simple rings with a minimal ideal,
{\em Ann. of Math.} (2) 90 (1969), 433--459.

\bibitem{Bl2} R.E. Block, On the extensions of Lie algebras,
{\em Canad. J. Math.} 20 (1968),
1439--1450.

\bibitem{G} W. de Graaf, Experimental Mathematics: http://www.expmath.org/expmath/volumes/14/14.1/deGraaf.pdf

\bibitem{Jac} N. Jacobson, Classes of restricted Lie algebras of characteristic $p$, II, {\em Duke M. J.} { 10} (1943), 107 - 121.

\bibitem{Sch} C. Schneider, http://arxiv.org/abs/math.RA/0406365

\bibitem{St04} H. Strade, Simple Lie Algebras over Fields of Positive Characteristic,  Expositions in Mathematics 38, {\em deGruyter }Berlin, (2004).

\bibitem{SF} H. Strade, R. Farnsteiner,
Modular Lie algebras and their representations.
Monographs and Textbooks in Pure and Applied Mathematics,
116. {\em Marcel Dekker, Inc.}, New York, (1988).

\bibitem{ZP} J. Patera, H. Zassenhaus, Solvable Lie Algebras of Dimension $\le 4$ Over Perfect Fields, Linear Algebra and its Applications 142 (1990), 1 - 17


\end{thebibliography}
\end{document}